\documentclass[12pt]{article}
\usepackage{amsmath,amsfonts,amssymb,amsthm}
\input amssym.def
\begin{document}
\def \Z{\Bbb Z}
\def \C{\Bbb C}
\def \R{\Bbb R}
\def \N{\Bbb N}

\def \A{{\mathcal{A}}}
\def \D{{\mathcal{D}}}
\def \E{{\mathcal{E}}}
\def \L{\mathcal{L}}
\def \S{{\mathcal{S}}}
\def \Q{\mathbf{Q}}
\def \wt{{\rm wt}}
\def \tr{{\rm tr}}
\def \span{{\rm span}}
\def \Res{{\rm Res}}
\def \Der{{\rm Der}}
\def \End{{\rm End}}
\def \Ind {{\rm Ind}}
\def \Irr {{\rm Irr}}
\def \Aut{{\rm Aut}}
\def \GL{{\rm GL}}
\def \Hom{{\rm Hom}}
\def \mod{{\rm mod}}
\def \ann{{\rm Ann}}
\def \ad{{\rm ad}}
\def \rank{{\rm rank}\;}
\def \<{\langle}
\def \>{\rangle}

\def \g{{\frak{g}}}
\def \h{{\hbar}}
\def \k{{\frak{k}}}
\def \sl{{\frak{sl}}}
\def \gl{{\frak{gl}}}

\def \be{\begin{equation}\label}
\def \ee{\end{equation}}
\def \bex{\begin{example}\label}
\def \eex{\end{example}}
\def \bl{\begin{lem}\label}
\def \el{\end{lem}}
\def \bt{\begin{thm}\label}
\def \et{\end{thm}}
\def \bp{\begin{prop}\label}
\def \ep{\end{prop}}
\def \br{\begin{rem}\label}
\def \er{\end{rem}}
\def \bc{\begin{coro}\label}
\def \ec{\end{coro}}
\def \bd{\begin{de}\label}
\def \ed{\end{de}}

\newcommand{\m}{\bf m}
\newcommand{\n}{\bf n}
\newcommand{\nno}{\nonumber}
\newcommand{\nord}{\mbox{\scriptsize ${\circ\atop\circ}$}}
\newtheorem{thm}{Theorem}[section]
\newtheorem{prop}[thm]{Proposition}
\newtheorem{coro}[thm]{Corollary}
\newtheorem{conj}[thm]{Conjecture}
\newtheorem{example}[thm]{Example}
\newtheorem{lem}[thm]{Lemma}
\newtheorem{rem}[thm]{Remark}
\newtheorem{de}[thm]{Definition}
\newtheorem{hy}[thm]{Hypothesis}
\makeatletter
\@addtoreset{equation}{section}
\def\theequation{\thesection.\arabic{equation}}
\makeatother
\makeatletter

\begin{center}
{\Large \bf Twisted modules for quantum vertex algebras }
\end{center}

\begin{center}
{Haisheng Li$^{*}$\footnote{Partially supported by NSF grant
DMS-0600189}, Shaobin Tan$^{**}$\footnote{Partially supported by a
China NSF grant 10671160}, and Qing
Wang$^{**}$\\
$\mbox{}^{*}$Department of Mathematical Sciences, Rutgers
University,
Camden, NJ 08102\\
$\mbox{}^{**}$Department of Mathematics, Xiamen University, Xiamen,
China}
\end{center}

\begin{abstract}
We study twisted modules for (weak) quantum vertex algebras and we
give a conceptual construction of (weak) quantum vertex algebras and
their twisted modules. As an application we construct and classify
irreducible twisted modules for a certain family of quantum vertex
algebras.
\end{abstract}

\section{Introduction}
In \cite{ek}, one of a series of papers on quantizations of Lie
bialgebras, Etingof and Kazhdan developed a fundamental theory of
quantum vertex operator algebras. Inspired by this theory, one of us
(H.Li) has extensively studied a notion of (weak) quantum vertex
algebra (see \cite{li-qva1}, \cite{li-qva2}, \cite{li-qva3}). While
quantum vertex operator algebras in the sense of \cite{ek} are
formal deformations of vertex algebras, weak quantum vertex algebras
are generalizations of vertex algebras and vertex super-algebras in
a certain direction. A weak quantum vertex algebra satisfies all the
axioms that define the notion of vertex algebra (cf. \cite{ll})
except the Jacobi identity axiom that is replaced by a braided
Jacobi identity axiom. An interesting feature of this theory is that
each weak quantum vertex algebra $V$ that is non-degenerate in the
sense of \cite{ek} gives rise to a canonical rational unitary
quantum Yang-Baxter operator with one spectral parameter on $V$.
While Jacobi identity amounts to (suitably defined) commutativity
and associativity, braided Jacobi identity amounts to associativity
and braided commutativity, called $\S$-locality. Thus weak quantum
vertex algebras are automatically nonlocal vertex algebras, namely
field algebras in the sense of \cite{bk}, which are vertex analogs
of noncommutative associative algebras.

While nonlocal vertex algebras are arguably too wild to be
interesting, weak quantum vertex algebras form a perfect class from
various points of view. (Note that in \cite{ek} quantum vertex
operator algebras are singled out {}from the wild family of braided
vertex operator algebras by using ``associativity.'') On the one
hand, this class is large enough to include those known
generalizations of vertex algebras such as vertex super-algebras and
vertex color-algebras and to include interesting new examples.  On
the other hand, it is tangible. Because of the associativity, weak
quantum vertex algebras admit a very nice representation theory,
just as ordinary vertex algebras do. It was proved in \cite{li-qva1}
that for any general vector space $W$, every what was called
$\S$-local subset of $\Hom (W,W((x)))$ canonically generates a weak
quantum vertex algebra with $W$ as a canonical faithful module.
(This generalizes the corresponding result of \cite{li-local}, which
states that any local subset generates a vertex algebra with $W$ as
a canonical faithful module.)

In the theory of vertex (operator) algebras, for a vertex algebra
$V$, in addition to (untwisted) modules, one also has so-called
twisted modules with respect to a finite-order automorphism of $V$.
(This is one of those new features of vertex algebras, compared with
classical Lie or associative algebras.) The notion of twisted module
was originated from the construction of the celebrated moonshine
module vertex operator algebra $V^{\natural}$ (see \cite{flm1},
\cite{flm}, \cite{le}) and it had played a critical role therein.
Twisted representations for general vertex operator algebras have
been extensively studied in literature (cf. \cite{li-twisted} and
\cite{dlm-twisted}). In particular, a conceptual construction of
vertex algebras and their twisted modules was given in
\cite{li-twisted}, generalizing the result of \cite{li-local}.

In this current paper, we study twisted modules for (weak) quantum
vertex algebras and we establish a conceptual construction of weak
quantum vertex algebras and their twisted modules by using twisted
vertex operators. The main result of this paper naturally
generalizes the corresponding results of \cite{li-twisted} and
\cite{li-qva1}. As an application we construct and classify
irreducible twisted modules for  quantum vertex algebras $V_{\bf Q}$
which were constructed in \cite{kl} from a square multiplicatively
skew complex matrix ${\bf Q}$.

To describe the general construction, let $W$ be a general vector
space and let $N$ be a positive integer. We consider the space
$(\End W)[[x^{1/N},x^{-1/N}]]$, which is naturally graded by the
abelian group $\Z/N\Z$. Denote by $\theta$ the associated linear
automorphism of order $N$ (defined by using the principal $N$-th
root of unity). Set
$$\E(W,N)=\oplus_{i=0}^{N-1}x^{\frac{i}{N}}\Hom(W,W((x))),$$
a graded subspace, on which $\theta$ acts. We say a subset $U$ of
$\E(W,N)$ is $\S$-local if for any $a(x),b(x)\in U$, there exist
$$c_{i}(x),d_{i}(c)\in U,\; f_{i}(x)\in \C((x)),\; i=1,\dots,r$$
such that
$$(x_{1}-x_{2})^{k}a(x_{1})b(x_{2})
=(x_{1}-x_{2})^{k}\sum_{i=1}^{r}f_{i}(x_{2}-x_{1})c_{i}(x_{2})d_{i}(x_{1})$$
for some nonnegative integer $k$. We prove that any $\S$-local
subset of $\E(W,N)$, which spans a graded subspace, canonically
generates a weak quantum vertex algebra on which $\theta$ acts as an
automorphism of period $N$, and that $W$ is a canonical
$\theta$-twisted module.

Let $\Q=(q_{ij})$ be an $r\times r$ complex matrix satisfying the
condition that $q_{ij}q_{ji}=1$ for $1\le i,j\le r$. Define
$\A_{\Q}$ to be the associative algebra with generators $X_{i,n},\;
Y_{i,n}$ for $1\le i\le r,\; n\in \Z$, subject to relations
\begin{eqnarray*}
& &X_{i,m}X_{j,n}=q_{ij}X_{j,n}X_{i,m},\ \ \ \
Y_{i,m}Y_{j,n}=q_{ij}Y_{j,n}Y_{i,m},\\
&&X_{i,m}Y_{j,n}-q_{ij}^{-1}Y_{j,n}X_{i,m}=\delta_{ij}\delta_{m+n+1,0}
\end{eqnarray*}
for $1\le i,j\le r,\; m,n\in \Z$. Associative algebras like
$\A_{\Q}$ have appeared in the study of noncommutative quantum field
theory. Notice that if $q_{ij}=1$ for $1\le i,j\le r$, $\A_{\Q}$,
which is the universal enveloping algebra of a Heisenberg Lie
algebra, is a Weyl algebra. If $q_{ij}=-1$ for $1\le i,j\le r$,
$\A_{\Q}$, which is the universal enveloping algebra of a Lie
superalgebra, is a Clifford algebra. Thus algebras $\A_{\Q}$ are
generalizations of Weyl algebras and Clifford algebras. Let $J$ be
the left ideal of $\A_{\bf Q}$ generated by $X_{i,m},\; Y_{i,n}$ for
$1\le i\le r,\; m,n\ge 0$. Set $V_{\bf Q}=\A_{\bf Q}/J$, a left
$\A_{\bf Q}$-module, and set ${\bf 1}=1+J\in V_{\bf Q}$. It was
proved in \cite{kl} (cf. \cite{li-qva2}) that there exists a
canonical quantum vertex algebra structure on $V_{\bf Q}$.

Let $N$ be any positive integer. Define an automorphism $\theta_{N}$
of $\A_{\Q}$ such that
$$\theta_{N}(X_{i,m})=\omega_{N}X_{i,m},\ \ \ \ \
\theta_{N}(Y_{i,m})=\omega_{N}^{-1}Y_{i,m} $$ for $1\le i\le r,\;
m\in \Z$, where $\omega_{N}$ denotes the principal primitive $N$-th
root of unity. This gives rise to an order-$N$ automorphism of the
quantum vertex algebra $V_{\Q}$. In this paper, we study
$\theta_{N}$-twisted $V_{\Q}$-modules. We construct
$\theta_{N}$-twisted $V_{\Q}$-modules by using a twisted algebra of
$\A_{\bf Q}$.  Since $V_{\Q}$ are generalizations of Clifford vertex
operator superalgebras, our results generalize the corresponding
results of \cite{ffr} (cf. \cite{xu}).

This paper is organized as follows: In Section 2, we define the
notion of $\sigma$-twisted module for a nonlocal vertex algebra $V$
with an automorphism $\sigma$ of finite order and we establish
several basic results. In Section 3, we give a general construction
of weak quantum vertex algebras and twisted modules. In Section 4,
we construct and classify irreducible twisted modules for certain
quantum vertex algebras of Zamolodchikov-Faddeev type.

\section{Twisted modules for nonlocal vertex algebras}
In this section we formulate and study a notion of twisted module
for a general nonlocal vertex algebra, and we establish certain
basic properties. For this paper, the scalar field is the field $\C$
of complex numbers and $\N$ denotes the set of nonnegative integers.

We begin with the notion of nonlocal vertex algebra (\cite{li-qva1};
cf. \cite{bk}, \cite{li-g1}):

\bd{dgva} {\em A {\em nonlocal vertex algebra} is a vector space $V$
equipped with a linear map
\begin{eqnarray}
Y(\cdot,x): V &\rightarrow & \Hom (V,V((x)))\subset (\End V)[[x,x^{-1}]]\nonumber\\
v&\mapsto& Y(v,x)=\sum_{n\in \Z}v_{n}x^{-n-1}\;\;\; (v_{n}\in \End
V),
\end{eqnarray}
and equipped with a vector ${\bf 1}\in V,$ satisfying the conditions
that for $v\in V,$
\begin{eqnarray}
& &Y({\bf 1},x)v=v,\\
& &Y(v,x){\bf 1}\in V[[x]]\;\;\mbox{ and }\;\; \lim_{x\rightarrow
0}Y(v,x){\bf 1}=v
\end{eqnarray}
and that for $u,v,w\in V$, there exists a nonnegative integer $l$
such that
\begin{eqnarray}
(x_{0}+x_{2})^{l}Y(u,x_{0}+x_{2})Y(v,x_{2})w=
(x_{0}+x_{2})^{l}Y(Y(u,x_{0})v,x_{2})w.
\end{eqnarray}}
\ed

The following was proved in \cite{li-g1}:

\bl{lcon1} Let $V$ be a nonlocal vertex algebra. Define
 a linear operator ${\cal{D}}$ on $V$ by
\begin{eqnarray}
{\cal{D}}(v)=v_{-2}{\bf 1} \;\;\;\mbox{ for }v\in V.
\end{eqnarray}
Then
\begin{eqnarray}\label{edproperty}
[{\cal{D}},Y(v,x)]=Y({\cal{D}}(v),x)={d\over dx}Y(v,x)
\;\;\;\mbox{ for }v\in V.
\end{eqnarray}
Furthermore, for $v\in V$,
\begin{eqnarray}
&
&e^{x{\cal{D}}}Y(v,x_{1})e^{-x{\cal{D}}}=Y(e^{xD}v,x_{1})=Y(v,x_{1}+x),
\label{econjugationformula1}\\
& &Y(v,x){\bf 1}=e^{x{\cal{D}}}v.\label{ecreationwithd}
\end{eqnarray}
\el

Let $U$ be a subset of a nonlocal vertex algebra $V$. Denote by
$\<U\>$ the {\em nonlocal vertex subalgebra of $V$ generated by
$U$}, which by definition is the smallest nonlocal vertex subalgebra
of $V$, containing $U$. From \cite{li-g1} we have:

\bl{lsubalgebragenerated} For any subset $U$ of $V$, the nonlocal
vertex subalgebra $\<U\>$ generated by $U$ is linearly spanned by
the vectors
\begin{eqnarray}\label{espannform}
u^{(1)}_{n_{1}}\cdots u^{(r)}_{n_{r}}{\bf 1}
\end{eqnarray}
for $r\ge 0,\; u^{(i)}\in U,\; n_{1},\dots,n_{r}\in \Z$. \el

The notion of weak quantum vertex algebra, defined in
\cite{li-qva1}, singles out an important class of nonlocal vertex
algebras.

\bd{dweak-qva} {\em A {\em weak quantum vertex algebra} is a vector
space $V$ equipped with a linear map
$$Y: V\rightarrow \Hom (V,V((x)));\;\; v\mapsto Y(v,x)$$
and equipped with a vector ${\bf 1}\in V$, satisfying the conditions
that for $v\in V$,
\begin{eqnarray}
& &Y({\bf 1},x)v=v,\\
& &Y(v,x){\bf 1}\in V[[x]]\;\;\mbox{ and } \; \lim_{x\rightarrow
0}Y(v,x){\bf 1}=v,
\end{eqnarray}
and that for $u,v\in V$, there exist $v^{(i)}, u^{(i)}\in V,
f_{i}(x)\in \C((x)), i=1,\ldots,r,$ such that
\begin{eqnarray}\label{e-Rjacobi}
& &x_{0}^{-1}\delta\left(\frac{x_{1}-x_{2}}{x_{0}}\right)
Y(u,x_{1})Y(v,x_{2})\nonumber\\
& &\ \ \ \hspace{1cm}
-x_{0}^{-1}\delta\left(\frac{x_{2}-x_{1}}{-x_{0}}\right)
\sum_{i=1}^{r} f_{i}(-x_{0})Y(v^{(i)},x_{2})Y(u^{(i)},x_{1})\nonumber\\
&=&x_{2}^{-1}\delta\left(\frac{x_{1}-x_{0}}{x_{2}}\right)
Y(Y(u,x_{0})v,x_{2})
\end{eqnarray}
(the {\em $\S$-Jacobi identity}).} \ed

Just as Jacobi identity does, $\S$-Jacobi identity implies weak
associativity, so that a weak quantum vertex algebra is
automatically a nonlocal vertex algebra. Furthermore, $\S$-Jacobi
identity also implies {\em $\S$-locality} (cf. \cite{ek}) in the
sense that for any $u,v\in V$, there exist $f_{i}(x)\in \C((x)),\;
u^{(i)},v^{(i)}\in V$, $i=1,\dots,r,$ such that
\begin{eqnarray}\label{eslocality}
(x_{1}-x_{2})^{k}Y(u,x_{1})Y(v,x_{2})
=(x_{1}-x_{2})^{k}\sum_{i=1}^{r}f_{i}(x_{2}-x_{1})
Y(v^{(i)},x_{2})Y(u^{(i)},x_{1})
\end{eqnarray}
for some nonnegative integer $k$. In fact (see \cite{li-qva1}),
$\S$-Jacobi identity is equivalent to weak associativity and
$\S$-locality, and hence a weak quantum vertex algebra is exactly a
nonlocal vertex algebra that satisfies $\S$-locality.

\br{rrecallqva} {\em Let $V$ be a nonlocal vertex algebra and let
$$u,v, u^{(1)},v^{(1)},\dots,u^{(r)},v^{(r)}\in V;
\; f_{1}(x),\dots, f_{r}(x)\in \C((x)). $$ It was proved in
\cite{li-qva1} (Corollary 5.3) that
\begin{eqnarray*}
(x_{1}-x_{2})^{k}Y(u,x_{1})Y(v,x_{2})
=(x_{1}-x_{2})^{k}\sum_{i=1}^{r}f_{i}(x_{2}-x_{1})Y(v^{(i)},x_{2})Y(u^{(i)},x_{1})
\end{eqnarray*}
holds for some $k\in \N$ if and only if
\begin{eqnarray}
Y(u,x)v=\sum_{i=1}^{r}f_{i}(-x)e^{x\D}Y(v^{(i)},-x)u^{(i)}.
\end{eqnarray}
Following Etingof-Kazhdan \cite{ek}, define a linear map
$$Y(x): \C((x))\otimes V\otimes V\rightarrow V((x))$$
by $$Y(x)(f(x)\otimes u\otimes v)=f(x)Y(u,x)v \ \ \ \ \mbox{ for
}f(x)\in \C((x)),\; u,v\in V.$$ Then a nonlocal vertex algebra $V$
is a weak quantum vertex algebra if and only if
\begin{eqnarray}\label{eskew-symmetry-2}
Y(u,x)v\in e^{x\D}Y(x)(\C((x))\otimes V\otimes V)
\end{eqnarray}
for all $u,v\in V$. } \er

In practice, the following technical result is very useful (cf.
\cite{li-qva2}, Lemma 2.7):

\bp{pwqva-generating} Let $V$ be a nonlocal vertex algebra. Suppose
that $U$ is a subspace of $V$ satisfying the conditions that
$V=\<U\>$ and that for any $u,v\in U$, there exist
$$\sum_{j=1}^{r}v^{(j)}\otimes u^{(j)}\otimes f_{j}(x)\in U\otimes
U\otimes \C((x))$$ and a nonnegative integer $k$ such that
\begin{eqnarray}
(x_{1}-x_{2})^{k}Y(u,x_{1})Y(v,x_{2})
=(x_{1}-x_{2})^{k}\sum_{j=1}^{r}f_{j}(x_{2}-x_{1})Y(v^{(j)},x_{2})Y(u^{(j)},x_{1}).
\end{eqnarray}
Then $V$ is a weak quantum vertex algebra. \ep

\begin{proof} In view of Remark \ref{rrecallqva}, we must prove that (\ref{eskew-symmetry-2})
holds for all $u,v\in V$. Let us introduce some technical notions.
For a subspace $A$ of $V$, we set
$$A^{(2)}=\mbox{span}\{ u_{n}v\;|\; u,v\in A,\; n\in \Z\}\subset V.$$
We say that an ordered pair $(A,B)$ of subspaces of $V$ is
$\S$-local if
$$Y(a,x)b\in e^{x\D}Y(x)(\C((x))\otimes B\otimes A)
\ \ \ \ \mbox{ for }a\in A,\; b\in B.$$ For
$F(x_{1},x_{2}),G(x_{1},x_{2})\in V[[x_{1}^{\pm 1},x_{2}^{\pm 1}]]$,
we define $F\sim G$ if
$$(x_{1}-x_{2})^{p}(x_{1}+x_{2})^{q}F=(x_{1}-x_{2})^{p}(x_{1}+x_{2})^{q}G$$
for some $p,q\in \N$. It is clear that the defined relation
``$\sim$'' is an equivalence relation.

We are going to prove that if an ordered pair $(A,P)$ of subspaces
of $V$ is $\S$-local, then $(A,P^{(2)})$ and $(A^{(2)},P)$ are
$\S$-local. Then it follows {}from this and induction that
$(\<U\>,\<U\>)$ is $\S$-local. Therefore, $V$ is a weak quantum
vertex algebra.

First, we prove that $(A,P^{(2)})$ is $\S$-local. Let $a\in A,\;
u,v\in P$. There exist $f_{i}(x), g_{ij}(x)\in \C((x)),\;
 a^{(i)}, a^{(ij)}\in A,\; u^{(i)},v^{(j)}\in P$ for $1\le i\le r,
1\le j\le s,$ such that
$$Y(a,x_{1})Y(u,x_{2})v\sim \sum_{i=1}^{r}f_{i}(x_{2}-x_{1})
Y(u^{(i)},x_{2})Y(a^{(i)},x_{1})v,$$
$$Y(a^{(i)},x)v=\sum_{j=1}^{s}g_{ij}(-x)e^{x\D}Y(v^{(j)},-x)a^{(ij)},$$
 and
\begin{eqnarray}
Y(u^{(i)},x_{2}-x_{1})Y(v^{(j)},-x_{1})a^{(ij)} \sim
Y(Y(u^{(i)},x_{2})v^{(j)},-x_{1})a^{(ij)}.
\end{eqnarray}
Using the $\D$-bracket-derivative property
(\ref{econjugationformula1}) and weak associativity we get
\begin{eqnarray*}\label{s-skew}
& &Y(a,x_{1})Y(u,x_{2})v\nonumber\\
&\sim &\sum_{i=1}^{r}f_{i}(x_{2}-x_{1})
Y(u^{(i)},x_{2})Y(a^{(i)},x_{1})v\nonumber\\
&\sim &\sum_{i=1}^{r}f_{i}(x_{2}-x_{1})
Y(u^{(i)},x_{2})\sum_{j=1}^{s}g_{ij}(-x_{1})e^{x_{1}\D}Y(v^{(j)},-x_{1})a^{(ij)}\nonumber\\
&\sim&\sum_{i=1}^{r}\sum_{j=1}^{s}f_{i}(x_{2}-x_{1})g_{ij}(-x_{1})
e^{x_{1}\D}Y(u^{(i)},x_{2}-x_{1})Y(v^{(j)},-x_{1})a^{(ij)}\nonumber\\
&\sim&\sum_{i=1}^{r}\sum_{j=1}^{s}f_{i}(-x_{1}+x_{2})g_{ij}(-x_{1})
e^{x_{1}\D}Y(Y(u^{(i)},x_{2})v^{(j)},-x_{1})a^{(ij)}.
\end{eqnarray*}
That is, there exists a nonnegative integer $k$ such that
\begin{eqnarray*}
&&(x_{1}^{2}-x_{2}^{2})^{k}Y(a,x_{1})Y(u,x_{2})v\nonumber\\
&=&(x_{1}^{2}-x_{2}^{2})^{k}\sum_{i=1}^{r}\sum_{j=1}^{s}f_{i}(-x_{1}+x_{2})g_{ij}(-x_{1})
e^{x_{1}\D}Y(Y(u^{(i)},x_{2})v^{(j)},-x_{1})a^{(ij)}.
\end{eqnarray*}
As both sides involve only finitely many negative powers of $x_{2}$,
multiplying both sides by $(x_{1}-x_{2})^{-k}(x_{1}+x_{2})^{-k}$, we
obtain
$$Y(a,x_{1})Y(u,x_{2})v=\sum_{i=1}^{r}\sum_{j=1}^{s}f_{i}(-x_{1}+x_{2})g_{ij}(-x_{1})
e^{x_{1}\D}Y(Y(u^{(i)},x_{2})v^{(j)},-x_{1})a^{(ij)}.$$ It follows
that $(A,P^{(2)})$ is $\S$-local.

Next, we prove that $(A^{(2)},P)$ is $\S$-local. Let $a,b\in A,\;
w\in P$. There exist $f_{i}(x), g_{ij}(x)\in \C((x)),\;
a^{(ij)},b^{(j)}\in A,\; w^{(i)}, w^{(ij)}\in P$ for $1\le i\le r,\;
1\le j\le s,$  such that
$$Y(b,x_{2})w=\sum_{i} e^{x\D}f_{i}(-x)Y(w^{(i)},-x)b^{(i)},$$
$$Y(a,x_{1})Y(w^{(i)},-x_{2})b^{(i)}\sim \sum_{j=1}^{s}g_{ij}(-x_{2}-x_{1})
Y(w^{(ij)},-x_{2})Y(a^{(ij)},x_{1})b^{(i)}.$$ Using weak
associativity we get
\begin{eqnarray*}\label{s-skew2}
& &Y(Y(a,x_{1})b,x_{2})w\nonumber\\
&\sim& Y(a,x_{1}+x_{2})Y(b,x_{2})w\nonumber\\
&\sim&Y(a,x_{1}+x_{2})\sum_{i=1}^{r}f_{i}(-x_{2})
e^{x_{2}\D}Y(w^{(i)},-x_{2})b^{(i)}\nonumber\\
&\sim&\sum_{i=1}^{r}f_{i}(-x_{2})e^{x_{2}\D}Y(a,x_{1})Y(w^{(i)},-x_{2})b^{(i)}\nonumber\\
&\sim&\sum_{i=1}^{r}\sum_{j=1}^{s}f_{i}(-x_{2})g_{ij}(-x_{2}-x_{1})e^{x_{2}\D}
Y(w^{(ij)},-x_{2})Y(a^{(ij)},x_{1})b^{(i)}.
\end{eqnarray*}
By a similar reasoning we obtain
$$Y(Y(a,x_{1})b,x_{2})w
=\sum_{i=1}^{r}\sum_{j=1}^{s}f_{i}(-x_{2})g_{ij}(-x_{2}-x_{1})e^{x_{2}\D}
Y(w^{(ij)},-x_{2})Y(a^{(ij)},x_{1})b^{(i)}.$$ This proves that
$(A^{(2)},P)$ is $\S$-local. Now, the proof is complete.
\end{proof}

Next, we formulate a notion of twisted module for a general nonlocal
vertex algebra. For two nonlocal vertex algebras $V$ and $K$, a {\em
homomorphism} of nonlocal vertex algebras from $V$ to $K$ is a
linear map $\sigma$ such that
$$\sigma({\bf 1})={\bf 1},\;\;\sigma (Y(u,x)v)=Y(\sigma(u),x)\sigma(v)
\ \ \ \ \mbox{ for } u,v\in V.$$ An {\em automorphism} of a nonlocal
vertex algebra $V$ is a bijective endomorphism of $V$ and an
automorphism of a weak quantum vertex algebra $V$ is an automorphism
of $V$ viewed as a nonlocal vertex algebra.

Let $V$ be a nonlocal vertex algebra and let $\sigma$ be an
automorphism of $V$ of {\em period $N$} in the sense that
$\sigma^{N}=1$. Then $\sigma$ acts semisimply on $V$ and
$$V=V^{0}\oplus V^{1} \oplus\cdots\oplus V^{N-1},$$
where $V^{k}$ is the eigenspace of $\sigma$ with eigenvalue
$\exp\left(\frac{2k\pi \sqrt{-1}}{N}\right)$. Set
$$\omega_{N}=\exp\left(\frac{2\pi \sqrt{-1}}{N}\right),$$
the principal primitive $N$-th root of unity.

\bd{dtwisted-module} {\em Let $V$ be a nonlocal vertex algebra and
let $\sigma$ be an automorphism of $V$ of period $N$. A {\em
$\sigma$-twisted $V$-module} is a vector space $W$ equipped with a
linear map
$$Y_{W}: V\rightarrow \Hom (W,W((x^{\frac{1}{N}})))\subset (\End W)[[x^{1/N},x^{-1/N}]],$$
satisfying the conditions that $Y_{W}({\bf 1},x)=1_{W}$ (the
identity operator on $W$),
\begin{eqnarray}\label{etwisted-inv}
Y_{W}(\sigma v,x)=\lim_{x^{\frac{1}{N}}\rightarrow
\omega_{N}^{-1}x^{\frac{1}{N}}}Y_{W}(v,x)\ \ \ \mbox{ for }v\in V,
\end{eqnarray}
and that for $u,v\in V$, there exists a nonnegative integer $k$ such
that
\begin{eqnarray}\label{ehalf-local}
(x_{1}-x_{2})^{k}Y_{W}(u,x_{1})Y_{W}(v,x_{2}) \in \Hom
(W,W((x_{1}^{\frac{1}{N}},x_{2}^{\frac{1}{N}}))),
\end{eqnarray}
 and
\begin{eqnarray}\label{etwisted-assoc-sub}
x_{0}^{k}Y_{W}(Y(u,x_{0})v,x_{2})=\left((x_{1}-x_{2})^{k}
Y_{W}(u,x_{1})Y_{W}(v,x_{2})\right)|_{x_{1}^{\frac{1}{N}}=(x_{2}+x_{0})^{\frac{1}{N}}}.
\end{eqnarray}
If $V$ is a weak quantum vertex algebra, we define a {\em
$\sigma$-twisted $V$-module} to be a $\sigma$-twisted module for $V$
viewed as a nonlocal vertex algebra.} \ed

We shall denote a twisted module $W$ by a pair $(W,Y_{W})$ whenever
it is necessary. Notice that the axiom (\ref{etwisted-inv}) amounts
to that
\begin{eqnarray}\label{eshift-mod}
x^{\frac{r}{N}}Y_{W}(u,x)\in (\End W)[[x,x^{-1}]]
\end{eqnarray}
for $u\in V^{r},\; 0\le r\le N-1$.

\br{rjustification} {\em We here explain certain substitutions
including the one which has appeared in the right hand side of
(\ref{etwisted-assoc-sub}). First, recall the standard formal
variable expansion convention
$$(x_{2}\pm x_{0})^{\alpha}=\sum_{i\ge
0}\binom{\alpha}{i}(\pm 1)^{i}x_{2}^{\alpha -i}x_{0}^{i}\in
x_{2}^{\alpha}\C[x_{2}^{-1}][[x_{0}]]$$ for $\alpha\in \C$. Then for
$f(x_{1},x_{2})=\sum_{m,n\in
k+\N}f(m,n)x_{1}^{\frac{m}{N}}x_{2}^{\frac{n}{N}}\in
W((x_{1}^{\frac{1}{N}},x_{2}^{\frac{1}{N}}))$ with $k\in \Z$, we
define
\begin{eqnarray*}
f(x_{1},x_{2})|_{x_{1}^{1/N}=(x_{2}\pm x_{0})^{1/N}} &=&\sum_{m,n\in
k+\N}f(m,n)(x_{2}\pm x_{0})^{\frac{m}{N}}x_{2}^{\frac{n}{N}}\\
&=&\sum_{m,n\in k+\N}\sum_{i\ge 0}\binom{\frac{m}{N}}{i}(\pm
1)^{i}f(m,n)x_{2}^{\frac{m+n}{N}-i}x_{0}^{i},
\end{eqnarray*}
which lies in $W((x_{2}^{\frac{1}{N}}))[[x_{0}]]$. Furthermore, we
 extend this substitution in the obvious way for $F(x_{1},x_{2})\in
\Hom (W,W((x_{1}^{\frac{1}{N}},x_{2}^{\frac{1}{N}})))$, where
$$F(x_{1},x_{2})|_{x_{1}^{1/N}=(x_{2}+x_{0})^{1/N}}
\in (\Hom (W,W((x_{2}^{\frac{1}{N}}))))[[x_{0}]].$$ One can show
that this substitution gives rise to an injective map. That is, for
$F(x_{1},x_{2}), G(x_{1},x_{2})\in \Hom
(W,W((x_{1}^{\frac{1}{N}},x_{2}^{\frac{1}{N}})))$, if
$$F(x_{1},x_{2})|_{x_{1}^{1/N}=(x_{2}+x_{0})^{1/N}}
=G(x_{1},x_{2})|_{x_{1}^{1/N}=(x_{2}+x_{0})^{1/N}},$$ then
$F(x_{1},x_{2})=G(x_{1},x_{2})$.

We shall need another substitution
\begin{eqnarray}
H(x_{1},x_{2})|_{x_{1}^{1/N}=(x_{0}+x_{2})^{1/N}} \in
W((x_{0}^{\frac{1}{N}}))((x_{2}^{\frac{1}{N}}))
\end{eqnarray}
for $H(x_{1},x_{2})\in
W((x_{1}^{\frac{1}{N}}))((x_{2}^{\frac{1}{N}}))$. One can show that
\begin{eqnarray}
\left(F(x_{1},x_{2})|_{x_{2}^{1/N}=(x_{1}-x_{0})^{1/N}}\right)|_{x_{1}^{1/N}=(x_{2}+x_{0})^{1/N}}
=F(x_{1},x_{2})|_{x_{1}^{1/N}=(x_{2}+x_{0})^{1/N}}
\end{eqnarray}
for $F(x_{1},x_{2})\in\Hom
(W,W((x_{1}^{\frac{1}{N}},x_{2}^{\frac{1}{N}})))$. } \er

Note that for $n\in \Z$, we have
$\frac{1}{N}\sum_{j=1}^{N}\omega_{N}^{jn}=0$ if $n\notin N\Z$ and
$\frac{1}{N}\sum_{j=1}^{N}\omega_{N}^{jn}=1$ if $n\in N\Z$. We have
(see \cite{flm})
\begin{eqnarray}
x_{1}^{-1}\delta\left(\frac{x_{2}+x_{0}}{x_{1}}\right)
=\frac{1}{N}\sum_{j=0}^{N-1}
x_{1}^{-1}\delta\left(\omega_{N}^{j}\left(\frac{x_{2}+x_{0}}{x_{1}}\right)^{\frac{1}{N}}\right).
\end{eqnarray}
Furthermore, for $f(x_{1},x_{2})\in \Hom
(W,W((x_{1}^{\frac{1}{N}},x_{2}^{\frac{1}{N}})))$, we have
\begin{eqnarray}
&&x_{1}^{-1}\delta\left(\frac{x_{2}+x_{0}}{x_{1}}\right)f(x_{1},x_{2})\nonumber\\
&=&\frac{1}{N}\sum_{j=0}^{N-1} x_{1}^{-1}\delta\left(\omega_{N}^{j}
\left(\frac{x_{2}+x_{0}}{x_{1}}\right)^{\frac{1}{N}}\right)f(x_{1},x_{2})\nonumber\\
&=&\frac{1}{N}\sum_{j=0}^{N-1} x_{1}^{-1}\delta\left(\omega_{N}^{j}
\left(\frac{x_{2}+x_{0}}{x_{1}}\right)^{\frac{1}{N}}\right)
f(x_{1},x_{2})|_{x_{1}^{1/N}=\omega_{N}^{j}(x_{2}+x_{0})^{1/N}}.
\end{eqnarray}

For simplicity, sometimes we shall simply use $x_{1}=x_{2}+x_{0}$
for the substitution
$x_{1}^{\frac{1}{N}}=(x_{2}+x_{0})^{\frac{1}{N}}$.

\bl{associatity} Let $V$ be a nonlocal vertex algebra and let
$\sigma$ be an automorphism of period $N$. In the presence of the
other axioms in the definition of a $\sigma$-twisted $V$-module, the
axioms (\ref{ehalf-local}) and (\ref{etwisted-assoc-sub}) can be
equivalently replaced with the property that for $u\in V^{r},\; v\in
V,\; w\in W$ with $0\le r\le N-1$, there exists $l\in \N$ such that
\begin{eqnarray}\label{etwisted-weakassoc}
(x_{2}+x_{0})^{l+\frac{r}{N}}Y_{W}(Y(u,x_{0})v,x_{2})w
=(x_{0}+x_{2})^{l+\frac{r}{N}} Y_{W}(u,x_{0}+x_{2})Y_{W}(v,x_{2})w.
\end{eqnarray}
 \el

\begin{proof} Assume that $(W,Y_{W})$ is a
$\sigma$-twisted $V$-module. For $u\in V^{r},\; v\in V,\; w\in W$
with $0\le r\le N-1$, with (\ref{eshift-mod}) and
(\ref{ehalf-local}), there exists $l\in \N$ such that
$$x_{1}^{l+\frac{r}{N}}(x_{1}-x_{2})^{k}Y_{W}(u,x_{1})Y_{W}(v,x_{2})w
\in W[[x_{1}]]((x_{2}^{\frac{1}{N}})).$$ Then using
(\ref{etwisted-assoc-sub}) we obtain
\begin{eqnarray*}
&&(x_{0}+x_{2})^{l+\frac{r}{N}}x_{0}^{k}Y_{W}(u,x_{0}+x_{2})Y_{W}(v,x_{2})w\nonumber\\
&&
=\left(x_{1}^{l+\frac{r}{N}}(x_{1}-x_{2})^{k}
Y_{W}(u,x_{1})Y_{W}(v,x_{2})w\right)|_{x_{1}=x_{0}+x_{2}}
\nonumber\\
&&=\left(x_{1}^{l+\frac{r}{N}}(x_{1}-x_{2})^{k}
Y_{W}(u,x_{1})Y_{W}(v,x_{2})w\right)|_{x_{1}=x_{2}+x_{0}}
\nonumber\\
&&=(x_{2}+x_{0})^{l+\frac{r}{N}}
\left((x_{1}-x_{2})^{k}Y_{W}(u,x_{1})Y_{W}(v,x_{2})w\right)|_{x_{1}=x_{2}+x_{0}}
\nonumber\\
&&=(x_{2}+x_{0})^{l+\frac{r}{N}}x_{0}^{k}Y_{W}(Y(u,x_{0})v,x_{2})w,
\nonumber
\end{eqnarray*}
{}from which (\ref{etwisted-weakassoc}) follows. On the other hand,
assume that (\ref{etwisted-weakassoc}) holds. Let $k\in \N$ be such
that $x^{k}Y(u,x)v\in V[[x]]$. Then the left hand side of
(\ref{etwisted-weakassoc}), after multiplied by $x_{0}^{k}$,
involves only nonnegative integer powers of $x_{0}$. Consequently,
we have
\begin{eqnarray*}
(x_{0}+x_{2})^{l+\frac{r}{N}}x_{0}^{k}Y_{W}(u,x_{0}+x_{2})Y_{W}(v,x_{2})w\in
W[[x_{0}]]((x_{2}^{\frac{1}{N}})).
\end{eqnarray*}
 Substituting
$x_{0}=x_{1}-x_{2}$ we get
\begin{eqnarray}\label{eproperty-xyz}
x_{1}^{l+\frac{r}{N}}(x_{1}-x_{2})^{k}Y_{W}(u,x_{1})Y_{W}(v,x_{2})w\in
W[[x_{1}]]((x_{2}^{\frac{1}{N}})).
\end{eqnarray}
 That is,
$$(x_{1}-x_{2})^{k}Y_{W}(u,x_{1})Y_{W}(v,x_{2})w\in
W((x_{1}^{\frac{1}{N}},x_{2}^{\frac{1}{N}})).$$ Noticing that $k$
does not depend on $w$, we obtain
$$(x_{1}-x_{2})^{k}Y_{W}(u,x_{1})Y_{W}(v,x_{2})\in
\Hom (W,W((x_{1}^{\frac{1}{N}},x_{2}^{\frac{1}{N}}))),$$ proving
(\ref{ehalf-local}). Also, using (\ref{eproperty-xyz}) and
(\ref{etwisted-weakassoc}) we get
\begin{eqnarray*}
&&(x_{2}+x_{0})^{l+\frac{r}{N}}x_{0}^{k}Y_{W}(Y(u,x_{0})v,x_{2})w\\
&=&(x_{0}+x_{2})^{l+\frac{r}{N}}x_{0}^{k}
Y_{W}(u,x_{0}+x_{2})Y_{W}(v,x_{2})w\\
&=&\left(x_{1}^{l+\frac{r}{N}}(x_{1}-x_{2})^{k}
Y_{W}(u,x_{1})Y_{W}(v,x_{2})w\right)|_{x_{1}=x_{0}+x_{2}}\\
&=&\left(x_{1}^{l+\frac{r}{N}}(x_{1}-x_{2})^{k}
Y_{W}(u,x_{1})Y_{W}(v,x_{2})w\right)|_{x_{1}=x_{2}+x_{0}}\\
&=&(x_{2}+x_{0})^{l+\frac{r}{N}}\left((x_{1}-x_{2})^{k}
Y_{W}(u,x_{1})Y_{W}(v,x_{2})w\right)|_{x_{1}=x_{2}+x_{0}}.
\end{eqnarray*}
Multiplying both sides by $(x_{2}+x_{0})^{-l-\frac{r}{N}}$ we obtain
(\ref{etwisted-assoc-sub}).
\end{proof}

Using an argument similar to that in the proof of Proposition 2.6 in
\cite{li-g1} (cf. \cite{dlm-twisted}) we get:

\bl{lquasi-dproperty} Let $(W,Y_{W})$ be a $\sigma$-twisted module
for a nonlocal vertex algebra $V$. Then
\begin{eqnarray}\label{eDproperty-twistedmodule}
Y_{W}(\D v,x)=\frac{d}{dx}Y_{W}(v,x)\;\;\;\mbox{ for }v\in V.
\end{eqnarray}
\el

Furthermore, we have:

\bl{S-Jacobi} Let $\sigma$ be an automorphism of a nonlocal vertex
algebra $V$ of period $N$ and let $(W,Y_{W})$ be a $\sigma$-twisted
$V$-module. Let
$$a,b,a^{(1)},b^{(1)},\dots,a^{(r)},b^{(r)}\in V;
\;\; f_{1}(x),\dots,f_{r}(x)\in \C((x)).$$ Then
 \begin{eqnarray}\label{emodule-comm-relation}
 &&(x_{1}-x_{2})^{k}Y_{W}(a,x_{1})Y_{W}(b,x_{2})\nonumber\\
&&=(x_{1}-x_{2})^{k}\sum_{i=1}^{r}f_{i}(x_{2}-x_{1})
Y_{W}(b^{(i)},x_{2})Y_{W}(a^{(i)},x_{1})
\end{eqnarray}
holds for some $k\in \N$ if and only if
\begin{eqnarray}\label{t-Jacobi-lemma}
& &x_{0}^{-1}\delta\left(\frac{x_{1}-x_{2}}{x_{0}}\right)
Y_{W}(a,x_{1})Y_{W}(b,x_{2})\nonumber\\
&&\ \ \ \ \ -x_{0}^{-1}\delta\left(\frac{x_{2}-x_{1}}{-x_{0}}\right)
\sum_{i=1}^{r}f_{i}(-x_{0}) Y_{W}(b^{(i)},x_{2})Y_{W}(a^{(i)},x_{1})
\nonumber\\
&=&\frac{1}{N}\sum_{j=0}^{N-1}
x_{2}^{-1}\delta\left(\omega_{N}^{-j}\left(\frac{x_{1}-x_{0}}{x_{2}}\right)^{\frac{1}{N}}\right)
Y_{W}(Y(\sigma^{j}a,x_{0})b,x_{2}).
\end{eqnarray}
\el

\begin{proof} Clearly, (\ref{emodule-comm-relation}) follows from (\ref{t-Jacobi-lemma}).
Now we prove that (\ref{emodule-comm-relation}) also implies
(\ref{t-Jacobi-lemma}).
 Choose a positive integer $k$ so large that (\ref{emodule-comm-relation})
holds and that
$$(x_{1}-x_{2})^{k}Y_{W}(\sigma^{j}a,x_{1})Y_{W}(b,x_{2})
\in \Hom (W,W((x_{1}^{\frac{1}{N}},x_{2}^{\frac{1}{N}})))$$ for
$0\le j\le N-1$. Then
$$x_{0}^{k}Y_{W}(Y(\sigma^{j}a,x_{0})b,x_{2})=\left(
(x_{1}-x_{2})^{k}Y_{W}(\sigma^{j}a,x_{1})Y_{W}(b,x_{2})\right)|_{x_{1}^{1/N}=(x_{2}+x_{0})^{1/N}}
$$
for $0\le j\le N-1$.  Since
$$Y_{W}(\sigma^{j}a,x_{1})=\lim_{x_{1}^{1/N}\rightarrow
\omega_{N}^{-j}x_{1}^{1/N}}Y_{W}(a,x_{1})$$ (recall
(\ref{eshift-mod})), we have
$$x_{0}^{k}Y_{W}(Y(\sigma^{j}a,x_{0})b,x_{2})=\left(
(x_{1}-x_{2})^{k}Y_{W}(a,x_{1})Y_{W}(b,x_{2})\right)|_{x_{1}^{1/N}=\omega_{N}^{-j}(x_{2}+x_{0})^{1/N}}.
$$
Then using (\ref{emodule-comm-relation}) and Remark
\ref{rjustification} we get
\begin{eqnarray*}
&&x_{0}^{k}x_{0}^{-1}\delta\left(\frac{x_{1}-x_{2}}{x_{0}}\right)
Y_{W}(a,x_{1})Y_{W}(b,x_{2})\nonumber\\
& &\ \ \ \
-x_{0}^{k}x_{0}^{-1}\delta\left(\frac{x_{2}-x_{1}}{-x_{0}}\right)
\sum_{i=1}^{r}f_{i}(-x_{0})
Y_{W}(b^{(i)},x_{2})Y_{W}(a^{(i)},x_{1})\nonumber\\
&=&x_{0}^{-1}\delta\left(\frac{x_{1}-x_{2}}{x_{0}}\right)
\left((x_{1}-x_{2})^{k}
Y_{W}(a,x_{1})Y_{W}(b,x_{2})\right)\nonumber\\
&&\ \ -x_{0}^{-1}\delta\left(\frac{x_{2}-x_{1}}{-x_{0}}\right)
\left((x_{1}-x_{2})^{k}\sum_{i=1}^{r}f_{i}(x_{2}-x_{1})
Y_{W}(b^{(i)},x_{2})Y_{W}(a^{(i)},x_{1})\right)\nonumber\\
&=&x_{1}^{-1}\delta\left(\frac{x_{2}+x_{0}}{x_{1}}\right)
\left((x_{1}-x_{2})^{k}
Y_{W}(a,x_{1})Y_{W}(b,x_{2})\right)\nonumber\\
&=&\frac{1}{N}\sum_{j=0}^{N-1}
x_{1}^{-1}\delta\left(\omega_{N}^{-j}\left(\frac{x_{2}+x_{0}}{x_{1}}\right)^{\frac{1}{N}}\right)
\\
&&\ \ \ \ \ \cdot \left((x_{1}-x_{2})^{k}
Y_{W}(a,x_{1})Y_{W}(b,x_{2})\right)|_{x_{1}^{1/N}=\omega_{N}^{-j}(x_{2}+x_{0})^{1/N}}\nonumber\\
&=&\frac{1}{N}\sum_{j=0}^{N-1}
x_{1}^{-1}\delta\left(\omega_{N}^{-j}\left(\frac{x_{2}+x_{0}}{x_{1}}\right)^{\frac{1}{N}}\right)
\\
&&\ \ \ \ \ \cdot \left((x_{1}-x_{2})^{k}
Y_{W}(\sigma^{j}a,x_{1})Y_{W}(b,x_{2})\right)|_{x_{1}^{1/N}=(x_{2}+x_{0})^{1/N}}\nonumber\\
&=&x_{0}^{k}\frac{1}{N}\sum_{j=0}^{N-1}
x_{1}^{-1}\delta\left(\omega_{N}^{-j}\left(\frac{x_{2}+x_{0}}{x_{1}}\right)^{\frac{1}{N}}\right)
Y_{W}(Y(\sigma^{j}a,x_{0})b,x_{2}).
\end{eqnarray*}
Multiplying both sides by $x_{0}^{-k}$, we obtain
(\ref{t-Jacobi-lemma}).\end{proof}

\bp{pslocal-twisted} Let $V$ be a nonlocal vertex algebra with an
automorphism $\sigma$ of period $N$ and let $(W,Y_{W})$ be a
$\sigma$-twisted $V$-module, and let
$$u,v, u^{(1)},v^{(1)},\dots,u^{(r)},v^{(r)}\in V,\;
f_{1}(x),\dots,f_{r}(x)\in \C((x)).$$  If
\begin{eqnarray}\label{eslocal-relation-a}
(x_{1}-x_{2})^{k}Y(u,x_{1})Y(v,x_{2})
=(x_{1}-x_{2})^{k}\sum_{i=1}^{r}f_{i}(x_{2}-x_{1})Y(v^{(i)},x_{2})Y(u^{(i)},x_{1})
\ \ \
\end{eqnarray}
for some $k\in \N$,  then
\begin{eqnarray}\label{eslocal-relation-m}
& &(x_{1}-x_{2})^{p}Y_{W}(u,x_{1})Y_{W}(v,x_{2})\nonumber\\
&=&(x_{1}-x_{2})^{p}\sum_{i=1}^{r}f_{i}(x_{2}-x_{1})Y_{W}(v^{(i)},x_{2})Y_{W}(u^{(i)},x_{1})
\end{eqnarray}
for some $p\in \N$. If $(W,Y_{W})$ is faithful, the converse is also
true.
 \ep

\begin{proof} First, note that from Remark \ref{rrecallqva},
(\ref{eslocal-relation-a}) is equivalent to
\begin{eqnarray}\label{eskew-proof}
Y(u,x_{0})v=\sum_{i=1}^{r}
f_{i}(-x_{0})e^{x_{0}\D}Y(v^{(i)},-x_{0})u^{(i)}.
\end{eqnarray}
Let $q\in \N$ be such that
\begin{eqnarray*}
& &(x_{1}-x_{2})^{q}Y_{W}(u,x_{1})Y_{W}(v,x_{2}) \in \Hom
(W,W((x_{1}^{1/N},x_{2}^{1/N}))),\\
&&(x_{1}-x_{2})^{q} Y_{W}(v^{(i)},x_{2})Y_{W}(u^{(i)},x_{1})\in \Hom
(W,W((x_{1}^{1/N},x_{2}^{1/N})))
\end{eqnarray*}
for all $1\le i\le r$. Furthermore, let $s\in \N$ be such that
$x^{s}f_{i}(x)\in \C[[x]]$ for $1\le i\le r$. Then
\begin{eqnarray}\label{associativity}
x_{0}^{q}Y_{W}(Y(u,x_{0})v,x_{2})=\left((x_{1}-x_{2})^{q}
Y_{W}(u,x_{1})Y_{W}(v,x_{2})\right)|_{x_{1}=x_{2}+x_{0}}
\end{eqnarray}
and
\begin{eqnarray}
&&x_{0}^{q+s}f_{i}(-x_{0})Y_{W}(Y(v^{(i)},-x_{0})u^{(i)},x_{1})\nonumber\\
&=&\left((x_{1}-x_{2})^{s}f_{i}(x_{2}-x_{1})\right)|_{x_{2}=x_{1}-x_{0}}
\left((x_{1}-x_{2})^{q}
Y_{W}(v^{(i)},x_{2})Y_{W}(u^{(i)},x_{1})\right)|_{x_{2}=x_{1}-x_{0}}\nonumber\\
&=&\left((x_{1}-x_{2})^{q+s}f_{i}(x_{2}-x_{1})
Y_{W}(v^{(i)},x_{2})Y_{W}(u^{(i)},x_{1})\right)|_{x_{2}=x_{1}-x_{0}}
\end{eqnarray}
for $1\le i\le r$. Then, using Lemma \ref{lquasi-dproperty} we
obtain
\begin{eqnarray}\label{eassociativity-other}
&&x_{0}^{q+s}\sum_{i=1}^{r}f_{i}(-x_{0})Y_{W}(e^{x_{0}\D}Y(v^{(i)},-x_{0})u^{(i)},x_{2})\nonumber\\
&=&x_{0}^{q+s}\sum_{i=1}^{r}f_{i}(-x_{0})Y_{W}(Y(v^{(i)},-x_{0})u^{(i)},x_{2}+x_{0})\nonumber\\
&=&\left(\sum_{i=1}^{r}(x_{1}-x_{2})^{q+s}f_{i}(x_{2}-x_{1})
Y_{W}(v^{(i)},x_{2})Y_{W}(u^{(i)},x_{1})\right)|_{x_{2}=x_{1}-x_{0}}|_{x_{1}=x_{2}+x_{0}}\nonumber\\
&=&\left(\sum_{i=1}^{r}(x_{1}-x_{2})^{q+s}f_{i}(x_{2}-x_{1})
Y_{W}(v^{(i)},x_{2})Y_{W}(u^{(i)},x_{1})\right)|_{x_{1}=x_{2}+x_{0}}.
\end{eqnarray}
If (\ref{eskew-proof}) holds, combining (\ref{associativity}) with
(\ref{eassociativity-other}) we obtain (\ref{eslocal-relation-m})
with $p=q+s$. On the other hand, if $(W,Y_{W})$ is faithful and if
(\ref{eslocal-relation-m}) holds, combining (\ref{associativity})
with (\ref{eassociativity-other}) we get (\ref{eskew-proof}).
\end{proof}

Furthermore, we have the following refinement of Proposition
\ref{pslocal-twisted}, which is a twisted analog of Proposition 6.7
of \cite{li-qva1}:

\bp{q-communication} Let $V$ be a nonlocal vertex algebra with an
automorphism $\sigma$ of period $N$ and let $(W,Y_{W})$ be a
$\sigma$-twisted $V$-module. Let
\begin{eqnarray*}
&& n\in\Z,\; u\in V^{k}, v, u^{(1)}, v^{(1)},\ldots, u^{(r)},
v^{(r)},\; c^{(0)},c^{(1)},\dots,c^{(s)}\in V,\\
&&\hspace{2cm} f_{1}(x),\ldots,f_{r}(x)\in \C((x)).
\end{eqnarray*}
If
\begin{eqnarray}\label{etranform-algebra}
&&(x_{1}-x_{2})^{n}Y(u,x_{1})Y(v,x_{2})
-(-x_{2}+x_{1})^{n}\sum_{i=1}^{r}f_{i}(x_{2}-x_{1})
Y(v^{(i)},x_{2})Y(u^{(i)},x_{1})\nonumber\\
&=&
\sum_{j=0}^{s}Y(c^{(j)},x_{2})\frac{1}{j!}\left(\frac{\partial}{\partial
x_{2}}\right)^{j}x_{2}^{-1}\delta\left(\frac{x_{1}}{x_{2}}\right),
\end{eqnarray}
then
\begin{eqnarray}\label{etranform-module}
&&(x_{1}-x_{2})^{n}Y_{W}(u,x_{1})Y_{W}(v,x_{2})\nonumber\\
&&\hspace{1cm}-(-x_{2}+x_{1})^{n}\sum_{i=1}^{r}f_{i}(x_{2}-x_{1})
Y_{W}(v^{(i)},x_{2})Y_{W}(u^{(i)},x_{1})\nonumber\\
&=&
\sum_{j=0}^{s}Y_{W}(c^{(j)},x_{2})\frac{1}{j!}\left(\left(\frac{\partial}{\partial
x_{2}}\right)^{j}x_{2}^{-1}\delta\left(\frac{x_{1}}{x_{2}}\right)
\left(\frac{x_{1}}{x_{2}}\right)^{\frac{k}{N}}\right).
\end{eqnarray}
If $(W,Y_{W})$ is faithful, the converse is also true.\ep

\begin{proof}  {}From the first part of the
proof of Proposition 6.7 of \cite{li-qva1}, we see that
(\ref{etranform-algebra}) is equivalent to
\begin{eqnarray}
&&\ \ \ \ Y(u,x)v=\sum_{i=1}^{r}f_{i}(-x)e^{x\D} Y(v^{(i)},-x)u^{(i)},
\label{eskew-proof2}\\
&&u_{n+j}v=c^{(j)}\ \ \ \mbox{ for }0\le j\le s \ \ \mbox{ and }\ \
u_{n+j}v=0\ \ \ \mbox{ for }\ j> s.\label{ecomponents}
\end{eqnarray}
Assume that (\ref{etranform-algebra}) holds. Then
(\ref{eskew-proof2}) and  (\ref{ecomponents}) hold. By Proposition
\ref{pslocal-twisted}, (\ref{eslocal-relation-m}) holds.
Furthermore, using Lemma \ref{S-Jacobi} we get
\begin{eqnarray}\label{ecommon}
&&(x_{1}-x_{2})^{n}Y_{W}(u,x_{1})Y_{W}(v,x_{2})\nonumber\\
&&\hspace{1cm}-(-x_{2}+x_{1})^{n}\sum_{i=1}^{r}f_{i}(x_{2}-x_{1})
Y_{W}(v^{(i)},x_{2})Y_{W}(u^{(i)},x_{1})\nonumber\\
&=&
\sum_{j=0}^{s}Y_{W}(u_{j+n}v,x_{2})\frac{1}{j!}\left(\left(\frac{\partial}{\partial
x_{2}}\right)^{j}x_{2}^{-1}\delta\left(\frac{x_{1}}{x_{2}}\right)
\left(\frac{x_{1}}{x_{2}}\right)^{\frac{k}{N}}\right).
\end{eqnarray}
Then using  (\ref{ecomponents}) we obtain (\ref{etranform-module}).
On the other hand, assume that $(W,Y_{W})$ is faithful and that
(\ref{etranform-module}) holds. It follows that
(\ref{eslocal-relation-m}) holds. From Lemma \ref{S-Jacobi}, we see
that (\ref{ecommon}) holds. With $W$ faithful, combining
(\ref{ecommon}) with (\ref{etranform-module}), we obtain
(\ref{ecomponents}). Also, by Proposition \ref{pslocal-twisted}, we
have (\ref{eslocal-relation-a}), which implies (\ref{eskew-proof2}).
Now, (\ref{etranform-algebra}) follows.
\end{proof}

As an immediate consequence we have:

\bc{cvirasoro} Let $V$ be a nonlocal vertex algebra with an
automorphism $\sigma$ of period $N$ and let $(W,Y_{W})$ be a
$\sigma$-twisted $V$-module. Suppose that $\omega\in V^{0}$ (i.e.,
$\sigma(\omega)=\omega$) is a conformal vector in the sense that
$$[L(m),L(n)]=(m-n)L(m+n)+\frac{1}{12}(m^{3}-m)\delta_{m+n,0}c$$
for $m,n\in \Z$, and $L(-1)=\D$, where $Y(\omega,x)=\sum_{n\in
\Z}L(n)x^{-n-2}$ and $c$ is a complex number. Then $W$ is  a module
for the Virasoro algebra with $L(m)$ acting as $L_{W}(m)$ for $m\in
\Z$ with the same central charge $c$, where
$Y_{W}(\omega,x)=\sum_{n\in \Z}L_{W}(n)x^{-n-2}$. \ec

 Combining Lemma \ref{S-Jacobi} with Proposition
\ref{pslocal-twisted} we immediately have:

\bc{cqva-twisted-module} Let $V$ be a weak quantum vertex algebra,
let $\sigma$ be an automorphism of period $N$ of $V$, and let
$(W,Y_{W})$ be a $\sigma$-twisted module for $V$ viewed as a
nonlocal vertex algebra. Let
$$u,v, u^{(1)},v^{(1)},\dots,u^{(r)},v^{(r)}\in V,\;
f_{1}(x),\dots,f_{r}(x)\in \C((x)).$$  If
\begin{eqnarray*}
(x_{1}-x_{2})^{k}Y(u,x_{1})Y(v,x_{2})
=(x_{1}-x_{2})^{k}\sum_{i=1}^{r}f_{i}(x_{2}-x_{1})Y(v^{(i)},x_{2})Y(u^{(i)},x_{1})
\ \ \
\end{eqnarray*}
for some $k\in \N$, then
\begin{eqnarray}\label{t-Jacobi2}
& &x_{0}^{-1}\delta\left(\frac{x_{1}-x_{2}}{x_{0}}\right)
Y_{W}(a,x_{1})Y_{W}(b,x_{2})\nonumber\\
& &\ \ \ \hspace{1cm}
-x_{0}^{-1}\delta\left(\frac{x_{2}-x_{1}}{-x_{0}}\right)
\sum_{i=1}^{r}f_{i}(-x_{0})
Y_{W}(b^{(i)},x_{2})Y_{W}(a^{(i)},x_{1})\ \ \ \
\nonumber\\
&=&\frac{1}{N}\sum_{j=0}^{N-1}
x_{2}^{-1}\delta\left(\omega_{N}^{-j}\left(\frac{x_{1}-x_{0}}{x_{2}}\right)^{\frac{1}{N}}\right)
Y_{W}(Y(\sigma^{j}a,x_{0})b,x_{2}).
\end{eqnarray}
\ec

\br{rtwisted-module-va} {\em  Let $V$ be a vertex algebra and let
$\sigma$ be an automorphism of $V$ of period $N$. Recall (see
\cite{flm}, \cite{ffr}, \cite{dong1}) that a {\em $\sigma$-twisted
$V$-module} is a vector space $W$ equipped with a linear map $Y_{W}:
V\rightarrow \Hom (W,W((x^{\frac{1}{N}})))$, satisfying the
conditions that $Y_{W}({\bf 1},x)=1_{W}$,
\begin{eqnarray}\label{ehomog-u}
x^{\frac{j}{N}}Y_{W}(u,x)\in \Hom (W,W((x)))\ \ \
\mbox{ for }u\in V^{j},\; 0\le j\le N-1,
\end{eqnarray}
 and that for $u,v\in V$,
\begin{eqnarray}\label{t-Jacobi}
& &x_{0}^{-1}\delta\left(\frac{x_{1}-x_{2}}{x_{0}}\right)
Y_{W}(u,x_{1})Y_{W}(v,x_{2})
-x_{0}^{-1}\delta\left(\frac{x_{2}-x_{1}}{-x_{0}}\right)Y_{W}(v,x_{2})Y_{W}(u,x_{1})
\nonumber\\
& &\ \ \  \ \ \ \ =\frac{1}{N}\sum_{j=0}^{N-1}
x_{2}^{-1}\delta\left(\omega_{N}^{-j}\left(\frac{x_{1}-x_{0}}{x_{2}}\right)^{\frac{1}{N}}\right)
Y_{W}(Y(\sigma^{j}u,x_{0})v,x_{2}).
\end{eqnarray}
In view of Corollary \ref{cqva-twisted-module}, a $\sigma$-twisted
module for $V$ viewed as a vertex algebra is the same as a
$\sigma$-twisted module for $V$ viewed as a nonlocal vertex
algebra.} \er

\bex{example-twisted} {\em Let $W$ be a vector space and let $N$ be
a positive integer. Set $V=(\End W)((x^{\frac{1}{N}}))$. The space
$V$ is naturally an associative algebra with identity $1$ and the
formal differential operator $\frac{d}{dx}$ is a derivation of $V$.
Then $V$ is a nonlocal vertex algebra with $1$ as the vacuum vector
and with
$$Y(a(x),x_{0})b(x)=\left( e^{x_{0}\frac{d}{dx}}a(x)\right)b(x)
=a(x+x_{0})b(x)\ \ \ \mbox{ for }a(x),b(x)\in V.$$ Define a linear
map $\theta: V\rightarrow V$ by
$$\theta (a(x))=\lim_{x^{\frac{1}{N}}\rightarrow
\omega_{N}^{-1}x^{\frac{1}{N}}}a(x)\ \ \ \ \mbox{ for }a(x)\in V.$$
It is clear that $\theta$ is an order-$N$ automorphism of $V$ viewed
as an associative algebra. It follows that $\theta$ is an order-$N$
automorphism of $V$ viewed as a nonlocal vertex algebra.
 Furthermore, $W$
is a $\theta$-twisted $V$-module with $Y_{W}(a(x),z)=a(z)$ for
$a(x)\in V$. Indeed, for $a(x),b(x)\in V$, we have
$$a(x_{1})b(x_{2})\in \Hom (W,W((x_{1}^{\frac{1}{N}},x_{2}^{\frac{1}{N}}))),$$ so that
$$Y_{W}(a(x),x_{1})Y_{W}(b(x),x_{2})=a(x_{1})b(x_{2})
\in \Hom (W,W((x_{1}^{\frac{1}{N}},x_{2}^{\frac{1}{N}}))),$$
\begin{eqnarray*}
Y_{W}(Y(a(x),x_{0})b(x),x_{2})=\left(Y(a(x),x_{0})b(x)\right)|_{x=x_{2}}
=\left(a(x_{1})b(x_{2})\right)|_{x_{1}=x_{2}+x_{0}}.
\end{eqnarray*}
 } \eex

\section{General construction of twisted modules}

In this section we shall give a conceptual construction of twisted
modules for weak quantum vertex algebras, generalizing the
construction in \cite{li-twisted} of twisted modules for vertex
superalgebras and the construction in \cite{li-qva1} of weak quantum
vertex algebras.

Let $W$ be a general vector space and let $N$ be a positive integer,
both of which are fixed throughout this section. Set
$$\Z_{N}=\Z/N\Z,$$
a cyclic abelian group of order $N$, and  for $j\in \Z$, set
$\bar{j}=j+N\Z\in \Z_{N}$. Consider the vector space $(\End
W)[[x^{\frac{1}{N}},x^{-\frac{1}{N}}]]$, which is naturally
$\Z_{N}$-graded with
$$(\End W)[[x^{\frac{1}{N}},x^{-\frac{1}{N}}]]
=\bigoplus_{\bar{j}\in \Z_{N}}x^{\frac{j}{N}}(\End W)[[x,x^{-1}]].$$
We define an order-$N$ linear automorphism $\theta$ of $(\End
W)[[x^{\frac{1}{N}},x^{-\frac{1}{N}}]]$ by
$$\theta
f(x)=\omega_{N}^{-j}f(x)\ \ \ \mbox{ for }f(x)\in
x^{\frac{j}{N}}(\End W)[[x,x^{-1}]] \ \mbox{with }j\in \Z.$$ For a
general (not necessarily homogeneous) $f(x)$, we have
$$\theta (f(x))=\lim_{x^{\frac{1}{N}}\rightarrow
\omega_{N}^{-1}x^{\frac{1}{N}}}f(x)\ \ \ \mbox{ for }f(x)\in (\End
W)[[x^{\frac{1}{N}},x^{-\frac{1}{N}}]].$$

Set
\begin{eqnarray} \E(W,N)=\Hom (W,W((x^{\frac{1}{N}}))),
\end{eqnarray}
a subspace of $(\End W)[[x^{\frac{1}{N}},x^{-\frac{1}{N}}]]$,
recalling that $$\E(W)=\Hom (W,W((x)))\subset (\End
W)[[x,x^{-1}]].$$ It is clear that $\E(W,N)$ is a graded subspace of
$(\End W)[[x^{\frac{1}{N}},x^{-\frac{1}{N}}]]$ with
\begin{eqnarray}
\E(W,N)=\bigoplus_{j=0}^{N-1} x^{-\frac{j}{N}}\E(W),
\end{eqnarray}
so that $\theta$ is also an order-$N$ linear automorphism of
$\E(W,N)$. For $0\leq j\leq N-1$, set
$$\E(W,N)^{j}=x^{-\frac{j}{N}}\E(W)=
\{f(x)\in\E(W,N)\;|\;\theta (f(x))=\omega_{N}^{j}f(x)\}.$$

\bd{dcompatibility} {\em An (ordered) finite sequence
$(a_{1}(x),\dots,a_{r}(x))$ in $\E(W,N)$ is said to be {\em
compatible} if there exists a nonnegative integer $k$ such that
\begin{eqnarray}
\left(\prod_{1\le i<j\le
r}(x_{i}-x_{j})^{k}\right)a_{1}(x_{1})\cdots a_{r}(x_{r}) \in \Hom
(W,W((x_{1}^{\frac{1}{N}},\dots,x_{r}^{\frac{1}{N}}))).
\end{eqnarray}
A subset  $U$ of $\E(W,N)$ is said to be {\em compatible} if every
finite sequence in $U$ is compatible.} \ed

In particular, an ordered pair $(a(x),b(x))$ in $\E(W,N)$ is
compatible if and only if there exists a nonnegative integer $k$
such that
\begin{eqnarray}\label{ecompatible-2}
(x_{1}-x_{2})^{k}a(x_{1})b(x_{2})\in \Hom
(W,W((x_{1}^{\frac{1}{N}},x_{2}^{\frac{1}{N}}))).
\end{eqnarray}

\bd{d-operation} {\em Let $(a(x),b(x))$ be a compatible (ordered)
pair in $\E(W,N)$. Define $a(x)_{n}b(x)\in (\End
W)[[x^{\frac{1}{N}},x^{-\frac{1}{N}}]]$ for $n\in \Z$ in terms of
the generating function
\begin{eqnarray}
Y_{\E}(a(x),x_{0})b(x)=\sum_{n\in \Z}(a(x)_{n}b(x)) x_{0}^{-n-1}
\end{eqnarray}
by
\begin{eqnarray}
& &Y_{\E}(a(x),x_{0})b(x) =x_{0}^{-k}
((x_{1}-x)^{k}a(x_{1})b(x))|_{x_{1}^{\frac{1}{N}}=
(x+x_{0})^{\frac{1}{N}}},\label{edef-3.15}
\end{eqnarray}
where $k$ is any nonnegative integer such that (\ref{ecompatible-2})
holds.}\ed

Just as before, for convenience we simply use $x_{1}= x+x_{0}$ for
the substitution $x_{1}^{\frac{1}{N}}= (x+x_{0})^{\frac{1}{N}}$. The
following two lemmas are straightforward consequences:

\bl{lexistence} Let $(a(x),b(x))$ be a compatible pair in $\E(W,N)$.
Then
\begin{eqnarray*}
a(x)_{n}b(x)\in \E (W,N)\;(=\Hom (W,W((x^{\frac{1}{N}}))))
\;\;\;\mbox{ for } n\in \Z.
\end{eqnarray*}
Furthermore, if $k$ is a nonnegative integer such that
\begin{eqnarray*}
(x_{1}-x_{2})^{k}a(x_{1})b(x_{2})\in
\Hom(W,W((x_{1}^{\frac{1}{N}},x_{2}^{\frac{1}{N}}))),
\end{eqnarray*}
then
\begin{eqnarray*}\label{e-truncation-ab}
a(x)_{n}b(x)=0\;\;\;\mbox{ for }n\ge k.
\end{eqnarray*}
\el

\bl{lvacuum} For any $a(x)\in \E(W,N)$, the sequences $(1_{W},
a(x))$ and $(a(x), 1_{W})$ are compatible and we have
\begin{eqnarray*}
& &Y_{\E}(1_{W},x_{0})a(x)=a(x),\label{e-vacuum-calculus}\\
& &Y_{\E}(a(x),x_{0})1_{W}=a(x+x_{0}) =e^{x_{0}{d\over
dx}}a(x).\label{e-creation-calculus}
\end{eqnarray*}
\el

We shall need the following technical result:

\bl{lproof-need} Let $(a_{i}(x),b_{i}(x))$ $(i=1,\dots,n)$ be
 compatible pairs in $\E(W,N)$. Suppose that
\begin{eqnarray}\label{esum-lemma}
\sum_{i=1}^{n}(z-x)^{k_{i}}a_{i}(z)b_{i}(x)\in \Hom
(W,W((z^{\frac{1}{N}},x^{\frac{1}{N}})))
\end{eqnarray}
for some nonnegative integers $k_{i}$. Then
\begin{eqnarray}\label{esum=sum}
\sum_{i=1}^{n}x_{0}^{k_{i}}Y_{\E}(a_{i}(x),x_{0})b_{i}(x)
=\left(\sum_{i=1}^{n}(z-x)^{k_{i}}a_{i}(z)b_{i}(x)\right)|_{z=x+x_{0}}.
\end{eqnarray}
\el

\begin{proof} Let $k$ be a nonnegative integer such that
$$(z-x)^{k}a_{i}(z)b_{i}(x)\in \Hom (W,W((z^{\frac{1}{N}},x^{\frac{1}{N}})))
\;\;\;\mbox{ for }1\le i\le n.$$
{}From Definition \ref{d-operation}, we have
$$x_{0}^{k}Y_{\E}(a_{i}(x),x_{0})b_{i}(x)
=\left((z-x)^{k}a_{i}(z)b_{i}(x)\right)\mid_{z=x+x_{0}}$$ for $1\le
i\le n$. Then using (\ref{esum-lemma}) we get
\begin{eqnarray}
&
&x_{0}^{k}\sum_{i=1}^{n}x_{0}^{k_{i}}Y_{\E}(a_{i}(x),x_{0})b_{i}(x)
\nonumber\\
&=&\sum_{i=1}^{n}x_{0}^{k_{i}}\left((z-x)^{k}a_{i}(z)b_{i}(x)\right)|_{z=x+x_{0}}
\nonumber\\
&=&\left((z-x)^{k}\sum_{i=1}^{n}(z-x)^{k_{i}}a_{i}(z)b_{i}(x)\right)|_{z=x+x_{0}}
\nonumber\\
&=&x_{0}^{k}\left(\sum_{i=1}^{n}(z-x)^{k_{i}}a_{i}(z)b_{i}(x)\right)|_{z=x+x_{0}},
\end{eqnarray}
which implies (\ref{esum=sum}) immediately.
\end{proof}

Let $U$ be a subspace of $\E(W)$ such that every ordered pair in
$U$ is compatible. We say $U$ is {\em closed} if
\begin{eqnarray}
a(x)_{n}b(x)\in U\;\;\;\mbox{ for }a(x),b(x)\in U,\; n\in \Z.
\end{eqnarray}
We are going to prove that any closed compatible subspace containing
$1_{W}$ of $\E(W,N)$ is a nonlocal vertex algebra. To achieve this
goal we first prove:

\bl{lclosed} Assume that $V$ is a subspace of $\E(W,N)$ such that
any sequence in $V$ of length $2$ or $3$ is compatible and such
that $V$ is closed. Let $\psi(x),\phi(x),\theta(x)\in V$ and let
$k$ be a nonnegative integer such that
\begin{eqnarray}
& &(x-y)^{k}\phi(x)\theta(y)\in \Hom
(W,W((x^{\frac{1}{N}},y^{\frac{1}{N}}))),
\label{e4.31}\\
& &(x-y)^{k}(x-z)^{k}(y-z)^{k} \psi(x)\phi(y)\theta(z)\in \Hom
(W,W((x^{\frac{1}{N}},y^{\frac{1}{N}},z^{\frac{1}{N}}))).\ \ \ \ \
\label{e4.32}
\end{eqnarray}
Then
\begin{eqnarray}
& &x_{1}^{k}x_{2}^{k}(x_{1}-x_{2})^{k}
Y_{\cal{E}}(\psi(x),x_{1})Y_{\cal{E}}(\phi(x),x_{2})\theta(x)\nonumber\\
&=&((y-x)^{k}(z-x)^{k}(y-z)^{k} \psi(y)\phi
(z)\theta(x))|_{y=x+x_{1},z=x+x_{2}}.
\end{eqnarray}
\el

\begin{proof} With (\ref{e4.31}), from Definition \ref{d-operation} we have
\begin{eqnarray}\label{ehphi-theta}
x_{2}^{k}Y_{\cal{E}}(\phi (x),x_{2})\theta (x)
=((z-x)^{k}\phi(z)\theta(x))|_{z=x+x_{2}},
\end{eqnarray}
which gives
\begin{eqnarray}\label{eright-left}
& &(y-x)^{k}(y-x-x_{2})^{k}x_{2}^{k}\psi(y)
Y_{\cal{E}}(\phi (x),x_{2})\theta(x)\nonumber\\
&=&((y-x)^{k}(y-z)^{k} (z-x)^{k}
\psi(y)\phi(z)\theta(x))|_{z=x+x_{2}}.
\end{eqnarray}
{}From (\ref{e4.32}) we see that the expression on the right-hand
side of (\ref{eright-left}) lies in  $$(\Hom (W,W((y^{\frac{1}{N}},
x^{\frac{1}{N}})))[[x_{2}]],$$ so does the expression on the
left-hand side. That is,
\begin{eqnarray*}
(y-x)^{k}(y-x-x_{2})^{k}x_{2}^{k}\psi(y)
Y_{\cal{E}}(\phi(x),x_{2})\theta(x) \in (\Hom
(W,W((y^{\frac{1}{N}},x^{\frac{1}{N}})))[[x_{2}]].
\end{eqnarray*}
In view of Lemma \ref{lproof-need}, by considering the coefficient
of each power of $x_{2}$, we have
\begin{eqnarray}
& &x_{1}^{k}(x_{1}-x_{2})^{k}x_{2}^{k}
Y_{\cal{E}}(\psi(x),x_{1})Y_{\cal{E}}(\phi(x),x_{2})\theta(x)
\nonumber\\
&=&((y-x)^{k}(y-x-x_{2})^{k}x_{2}^{k}
\psi(y)Y_{\cal{E}}(\phi(x),x_{2})\theta(x))|_{y=x+x_{1}}.
\end{eqnarray}
Using this and (\ref{ehphi-theta}) we obtain
\begin{eqnarray*}
& &x_{2}^{k}x_{1}^{k}(x_{1}-x_{2})^{k}
Y_{\cal{E}}(\psi(x),x_{1})Y_{\cal{E}}(\phi(x),x_{2})\theta(x)\nonumber\\
&=&x_{2}^{k}(y-x)^{k}(y-x-x_{2})^{k}
\psi(y)Y_{\cal{E}}(\phi (x),x_{2})\theta(x))|_{y=x+x_{1}}\nonumber\\
&=&((z-x)^{k}(y-x)^{k}(y-z)^{k} \psi(y)\phi
(z)\theta(x))|_{y=x+x_{1},z=x+x_{2}},
\end{eqnarray*}
concluding the proof.
\end{proof}

For any formal series $a(x)=\sum_{n\in \Z}a_{n}x^{-n-1}$ (with
coefficients $a_{n}$ in any vector space) and for any $m\in \Z$, we
set
\begin{eqnarray}
a(x)_{\ge m}=\sum_{n\ge m}a_{n}x^{-n-1}.
\end{eqnarray}
Then for any polynomial $q(x)$ we have
\begin{eqnarray}\label{esimplefactresidule}
\Res_{x}x^{m}q(x)a(x)=\Res_{x}x^{m}q(x)a(x)_{\ge m}.
\end{eqnarray}

Now we are in a position to present our first key result:

\bt{tclosed} Let $V$ be a $\Z_{N}$-graded subspace of $\E(W,N)$ such
that any sequence in $V$ of length $2$ or $3$ is compatible and such
that $V$ contains $1_{W}$ and is closed. Then $(V,Y_{\E},1_{W})$
carries the structure of a nonlocal vertex algebra with $\theta$ as
an automorphism of period $N$, and $W$ is a faithful
$\theta$-twisted $V$-module with $Y_{W}(a(x),x_{0})=a(x_{0})$ for
$a(x)\in V$. \et

\begin{proof} For the assertion on the nonlocal vertex algebra structure,
with Lemmas \ref{lexistence} and \ref{lvacuum}, it remains to prove
weak associativity, i.e., for $\psi(x),\phi(x),\theta(x)\in V$,
there exists a nonnegative integer $k$ such that
\begin{eqnarray*}\label{eweakassocmainthem}
&&(x_{0}+x_{2})^{k}Y_{\cal{E}}(\psi(x),x_{0}+x_{2})Y_{\cal{E}}(\phi(x),x_{2})\theta(x)\\
&=&(x_{0}+x_{2})^{k}Y_{\cal{E}}(Y_{\cal{E}}(\psi(x),x_{0})\phi(x),x_{2})\theta(x).
\end{eqnarray*}
Let $k$ be a nonnegative integer such that
\begin{eqnarray*}
& &(x-y)^{k}\psi(x)\phi(y)\in \Hom (W,W((x^{\frac{1}{N}},y^{\frac{1}{N}}))),\\
& &(x-y)^{k}\phi(x)\theta(y)\in \Hom (W,W((x^{\frac{1}{N}},y^{\frac{1}{N}}))),\\
& &(x-y)^{k}(x-z)^{k}(y-z)^{k} \psi(x)\phi(y)\theta(z)\in \Hom
(W,W((x^{\frac{1}{N}},y^{\frac{1}{N}},z^{\frac{1}{N}}))).
\end{eqnarray*}
By Lemma \ref{lclosed}, we have
\begin{eqnarray}\label{e5.76}
& &x_{2}^{k}(x_{0}+x_{2})^{k}x_{0}^{k}
Y_{\cal{E}}(\psi(x),x_{0}+x_{2})Y_{\cal{E}}(\phi(x),x_{2})\theta(x)
\nonumber\\
&=&((z-x)^{k}(y-x)^{k}(y-z)^{k}
\psi(y)\phi(z)\theta(x))|_{y=x+x_{0}+x_{2},z=x+x_{2}}.
\end{eqnarray}

On the other hand, let $n\in \Z$ be {\em arbitrarily fixed}. Since
$\psi(x)_{m}\phi(x)=0$ for $m$ sufficiently large, there exists a
nonnegative integer $l$, depending on $n$, such that
\begin{eqnarray}\label{eges}
x_{2}^{l}(Y_{\cal{E}}(\psi(x)_{m}\phi(x), x_{2})\theta(x)
=((z-x)^{l}(\psi(z)_{m}\phi(z))\theta(x))|_{z=x+x_{2}}
\end{eqnarray}
for {\em all} $m\ge n$. With $(x-y)^{k}\psi(x)\phi(y)\in \Hom
(W,W((x^{\frac{1}{N}},y^{\frac{1}{N}})))$, from Definition
\ref{d-operation} we have
\begin{eqnarray}\label{ethis}
x_{0}^{k}(Y_{\cal{E}}(\psi(x_{2}),x_{0})\phi(x_{2}))\theta(x)
=((y-x_{2})^{k}\psi(y)\phi(x_{2})\theta(x))|_{y=x_{2}+x_{0}}.
\end{eqnarray}
Using (\ref{esimplefactresidule}), (\ref{eges}) and (\ref{ethis})
we get
\begin{eqnarray}\label{e5.78}
& &\Res_{x_{0}}x_{0}^{n}(x_{0}+x_{2})^{k}x_{0}^{k}
x_{2}^{l}Y_{\cal{E}}(Y_{\cal{E}}(\psi(x),x_{0})\phi(x), x_{2})\theta(x)\nonumber\\
&=&\Res_{x_{0}}x_{0}^{n}(x_{0}+x_{2})^{k} x_{0}^{k}x_{2}^{l}
Y_{\cal{E}}(Y_{\cal{E}}(\psi(x),x_{0})_{\ge n}\phi(x), x_{2})\theta(x)\nonumber\\
&=&\Res_{x_{0}}x_{0}^{n} (x_{0}+x_{2})^{k}x_{0}^{k}((z-x)^{l}
(Y_{\cal{E}}(\psi(z),x_{0})_{\ge n}\phi(z))\theta(x))|_{z=x+x_{2}}\nonumber\\
&=&\Res_{x_{0}}x_{0}^{n} (x_{0}+x_{2})^{k}x_{0}^{k}
((z-x)^{l}Y_{\cal{E}}(\psi(z),x_{0})\phi(z))\theta(x))|_{z=x+x_{2}}\nonumber\\
&=&\Res_{x_{0}}x_{0}^{n}((z+x_{0}-x)^{k}x_{0}^{k}
(z-x)^{l}(Y_{\cal{E}}(\psi(z),x_{0})\phi(z))\theta(x))|_{z=x+x_{2}}\nonumber\\
&=&\Res_{x_{0}}x_{0}^{n}((y-x)^{k}(y-z)^{k}
(z-x)^{l}\psi(y)\phi(z)\theta(x))|_{y=z+x_{0},z=x+x_{2}}
\nonumber\\
&=&\Res_{x_{0}}x_{0}^{n}((y-x)^{k}(y-z)^{k}
(z-x)^{l}\psi(y)\phi(z)\theta(x))|_{y=x+x_{2}+x_{0},z=x+x_{2}}.
\end{eqnarray}
Combining (\ref{e5.78}) with (\ref{e5.76}) we get
\begin{eqnarray}\label{e5.79}
& &\Res_{x_{0}}x_{0}^{n} x_{2}^{k}(x_{0}+x_{2})^{k}x_{0}^{k}
 x_{2}^{l}Y_{\cal{E}}(\psi(x),x_{0}+x_{2})
Y_{\cal{E}}(\phi(x),x_{2})\theta (x)\nonumber\\
&=&\Res_{x_{0}}x_{0}^{n}x_{2}^{k}(x_{0}+x_{2})^{k}x_{0}^{k}
x_{2}^{l} Y_{\cal{E}}(Y_{\cal{E}}(\psi(x),x_{0})\phi(x),
x_{2})\theta(x).
\end{eqnarray}
That is
\begin{eqnarray*}
& &\Res_{x_{0}}x_{0}^{n+k}(x_{0}+x_{2})^{k}
Y_{\cal{E}}(\psi(x),x_{0}+x_{2})Y_{\cal{E}}(\phi(x),x_{2})\theta(x)\nonumber\\
&=&\Res_{x_{0}}x_{0}^{n+k}(x_{0}+x_{2})^{k}
Y_{\cal{E}}(Y_{\cal{E}}(\psi(x),x_{0})\phi(x), x_{2})\theta(x).\ \
\ \
\end{eqnarray*}
Since  $n$ is arbitrary, we obtain
\begin{eqnarray*}\label{enearfinal-new}
& &(x_{0}+x_{2})^{k}
Y_{\cal{E}}(\psi(x),x_{0}+x_{2})Y_{\cal{E}}(\phi(x),x_{2})\theta(x)\nonumber\\
&=&(x_{0}+x_{2})^{k} Y_{\cal{E}}(Y_{\cal{E}}(\psi(x),x_{0})\phi(x),
x_{2})\theta(x),
\end{eqnarray*}
as desired. Thus $(V,Y_{\E},1_{W})$ carries the structure of a
nonlocal vertex algebra.

To show that $\theta$ is an automorphism of $V$, first we see that
$\theta(1_{W})=1_{W}$. Next, for $a(x),b(x)\in V$,  there exists
$k\in \N$ such that
\begin{eqnarray*}
&&(x_{1}-x_{2})^{k}a(x_{1})b(x_{2})\in \Hom
(W,W((x_{1}^{\frac{1}{N}},x_{2}^{\frac{1}{N}}))),\\
&&(x_{1}-x_{2})^{k}\theta(a(x_{1}))\theta(b(x_{2}))\in \Hom
(W,W((x_{1}^{\frac{1}{N}},x_{2}^{\frac{1}{N}}))).
\end{eqnarray*}
{}From Definition \ref{d-operation} we have
\begin{eqnarray*}
&&x_{0}^{k}Y_{\E}(a(x),x_{0})b(x)
=\left((x_{1}-x)^{k}a(x_{1})b(x)\right)|_{x_{1}^{1/N}=(x+x_{0})^{1/N}},\\
&&x_{0}^{k}Y_{\E}(\theta(a(x)),x_{0})\theta(b(x))
=\left((x_{1}-x)^{k}\theta(a(x_{1}))\theta(b(x))\right)|_{x_{1}^{1/N}=(x+x_{0})^{1/N}}.
\end{eqnarray*}
Then
\begin{eqnarray*}
&&x_{0}^{k}\theta\left(Y_{\E}(a(x),x_{0})b(x)\right)\\
&=&\lim_{x^{1/N}\rightarrow \omega_{N}^{-1}x^{1/N}}\left(x_{0}^{k}Y_{\E}(a(x),x_{0})b(x)\right)\\
&=&\lim_{x^{1/N}\rightarrow
\omega_{N}^{-1}x^{1/N}}\left(\left((x_{1}-x)^{k}a(x_{1})b(x)\right)|_{x_{1}^{1/N}=(x+x_{0})^{1/N}}\right)\\
&=&\left((x_{1}-x)^{k}\theta(a(x_{1}))\theta(b(x))\right)|_{x_{1}^{1/N}=(x+x_{0})^{1/N}}\\
&=&x_{0}^{k}Y_{\E}(\theta(a(x)),x_{0})\theta(b(x)),
\end{eqnarray*}
which implies
$$\theta\left(Y_{\E}(a(x),x_{0})b(x)\right)
=Y_{\E}(\theta(a(x)),x_{0})\theta(b(x)).$$ Thus $\theta$ is an
automorphism of $V$.

By definition we have $Y_{W}(1_{W},x)=1_{W}$ and
$$Y_{W}(\theta(a(x)),z)=\lim_{z^{1/N}\rightarrow \omega_{N}^{-1}z^{1/N}}a(z)
=\lim_{z^{1/N}\rightarrow \omega_{N}^{-1}z^{1/N}}Y_{W}(a(x),z)$$ for
$a(x)\in V$. Furthermore, for $a(x),b(x)\in V$, there exists a
nonnegative integer $k$ such that
$$(x_{1}-x_{2})^{k}a(x_{1})b(x_{2})\in \Hom
(W,W((x_{1}^{1/N},x_{2}^{1/N})))$$  and
$$x_{0}^{k}Y_{\E}(a(x),x_{0})b(x)
=\left((x_{1}-x)^{k}a(x_{1})b(x)\right)|_{x_{1}^{1/N}=(x+x_{0})^{1/N}}.$$
Then
\begin{eqnarray*}
& &x_{0}^{k}Y_{W}(Y_{\E}(a(x),x_{0})b(x),x_{2})\nonumber\\
&=&x_{0}^{k}\left(Y_{\E}(a(x),x_{0})b(x)\right)|_{x=x_{2}}\nonumber\\
&=&\left((x_{1}-x_{2})^{k}a(x_{1})b(x_{2})\right)|_{x_{1}^{1/N}=(x_{2}+x_{0})^{1/N}}\nonumber\\
 &=&
\left((x_{1}-x_{2})^{k}Y_{W}(a(x),x_{1})Y_{W}(b(x),x_{2})\right)|_{x_{1}^{1/N}=
(x_{2}+x_{0})^{1/N}}.
\end{eqnarray*}
This proves that $W$ is a $\theta$-twisted $V$-module with
$Y_{W}(a(x),x_{0})=a(x_{0})$ for $a(x)\in V.$ It is clear that $W$
is faithful.
\end{proof}

Next we are going to prove that any compatible subset of $\E(W,N)$
generates a nonlocal vertex algebra. The following is a key result:

\bp{pgeneratingcomplicatedone} Let $\psi_{1}(x),\dots,\psi_{r}(x),
a(x),b(x),\phi_{1}(x), \dots,\phi_{s}(x) \in \E(W,N)$. Assume that
the ordered sequences $(a(x), b(x))$ and
$$(\psi_{1}(x),\dots,\psi_{r}(x), a(x),b(x),\phi_{1}(x),
\dots,\phi_{s}(x))$$
 are compatible.
Then for any $n\in \Z$, the ordered sequence
$$(\psi_{1}(x),\dots,\psi_{r}(x),a(x)_{n}b(x),\phi_{1}(x),\dots,\phi_{s}(x))$$
 is compatible.
\ep

\begin{proof}
Let $k$ be a nonnegative integer such that
$$(x_{1}-x_{2})^{k}a(x_{1})b(x_{2})\in \Hom (W,W((x_{1}^{\frac{1}{N}},x_{2}^{\frac{1}{N}})))$$
and
\begin{eqnarray}\label{elong-exp}
& &\left(\prod_{1\le i<j\le r}(y_{i}-y_{j})^{k}\right)
\left(\prod_{1\le i\le r, 1\le j\le s}(y_{i}-z_{j})^{k}\right)
\left(\prod_{1\le i<j\le s}(z_{i}-z_{j})^{k}\right)\nonumber\\
& &\;\;\cdot (x_{1}-x_{2})^{k}
\left(\prod_{i=1}^{r}(x_{1}-y_{i})^{k}(x_{2}-y_{i})^{k}\right)
\left(\prod_{i=1}^{s}(x_{1}-z_{i})^{k}(x_{2}-z_{i})^{k}\right)\nonumber\\
& &\;\;\cdot \psi_{1}(y_{1})\cdots \psi_{r}(y_{r})
a(x_{1})b(x_{2})\phi_{1}(z_{1})\cdots \phi_{s}(z_{s})\nonumber\\
& &\in \Hom (W,W((y_{1}^{\frac{1}{N}},\dots,
y_{r}^{\frac{1}{N}},x_{1}^{\frac{1}{N}},x_{2}^{\frac{1}{N}},z_{1}^{\frac{1}{N}},\dots,z_{s}^{\frac{1}{N}}))).
\end{eqnarray}
Set
$$P=\prod_{1\le i<j\le r}(y_{i}-y_{j})^{k},\;\;\;\;
Q=\prod_{1\le i<j\le s}(z_{i}-z_{j})^{k},\;\;\;\; R=\prod_{1\le
i\le r,\; 1\le j\le s}(y_{i}-z_{j})^{^{k}}.$$ Let $n\in \Z$ be
{\em arbitrarily fixed}. There exists a nonnegative integer $l$
such that
\begin{eqnarray}\label{etruncationpsiphi}
x_{0}^{l+n-k} \in \C[[x_{0}]].
\end{eqnarray}
Using (\ref{etruncationpsiphi}) and Definition \ref{d-operation}
we obtain
\begin{eqnarray}\label{ecompatibilitythreeproof}
& &\prod_{i=1}^{r}(x_{2}-y_{i})^{kl}
\prod_{j=1}^{s}(x_{2}-z_{j})^{kl}\cdot
\psi_{1}(y_{1})\cdots\psi_{r}(y_{r})
(a(x_{2})_{n}b(x_{2}))\phi_{1}(z_{1})\cdots\phi_{s}(z_{s})\nonumber\\
&=&\Res_{x_{0}}x_{0}^{n}\prod_{i=1}^{r}(x_{2}-y_{i})^{kl}
\prod_{j=1}^{s}(x_{2}-z_{j})^{kl} \cdot\nonumber\\
&&\ \ \ \ \cdot \psi_{1}(y_{1})\cdots \psi_{r}(y_{r})
(Y_{\cal{E}}(a(x_{2}),x_{0})b(x_{2}))\phi_{1}(z_{1})\cdots
\phi_{s}(z_{s})
\nonumber\\
&=&\Res_{x_{0}}x_{0}^{n}\prod_{i=1}^{r}(x_{2}-y_{i})^{kl}
\prod_{j=1}^{s}(x_{2}-z_{j})^{kl}\nonumber\\
& &\cdot( x_{0}^{-k}(x_{1}-x_{2})^{k} \psi_{1}(y_{1})\cdots
\psi_{r}(y_{r}) a(x_{1})b(x_{2})\phi_{1}(z_{1})\cdots
\phi_{s}(z_{s}))|_{x_{1}=x_{2}+x_{0}}\nonumber\\
&=&\Res_{x_{0}}x_{0}^{n}(\prod_{i=1}^{r}(x_{1}-x_{0}-y_{i})^{kl}
\prod_{j=1}^{s}(x_{1}-x_{0}-z_{j})^{kl}\nonumber\\
& &\cdot x_{0}^{-k}(x_{1}-x_{2})^{k} \psi_{1}(y_{1})\cdots
\psi_{r}(y_{r}) a(x_{1})b(x_{2})\phi_{1}(z_{1})\cdots
\phi_{s}(z_{s}))|_{x_{1}=x_{2}+x_{0}}\nonumber\\
&=&\Res_{x_{0}}x_{0}^{n}( e^{-x_{0}\frac{\partial}{\partial
x_{1}}} \left(\prod_{i=1}^{r}(x_{1}-y_{i})^{k}
\prod_{j=1}^{s}(x_{1}-z_{j})^{k}\right)^{l}\nonumber\\
&&\cdot x_{0}^{-k}(x_{1}-x_{2})^{k}\psi_{1}(y_{1})\cdots
\psi_{r}(y_{r}) a(x_{1})b(x_{2})\phi_{1}(z_{1})\cdots
\phi_{s}(z_{s}))|_{x_{1}=x_{2}+x_{0}}\nonumber\\
&=&\Res_{x_{0}}(\sum_{t=0}^{l-1}\frac{(-1)^{t}}{t!}x_{0}^{n+t}
\left(\frac{\partial}{\partial x_{1}}\right)^{t}
\left(\prod_{i=1}^{r}(x_{1}-y_{i})^{k}
\prod_{j=1}^{s}(x_{1}-z_{j})^{k}\right)^{l}\nonumber\\
& &\cdot x_{0}^{-k}(x_{1}-x_{2})^{k} \psi_{1}(y_{1})\cdots
\psi_{r}(y_{r}) a(x_{1})b(x_{2})\phi_{1}(z_{1})\cdots
\phi_{s}(z_{s}))|_{x_{1}=x_{2}+x_{0}}.
\end{eqnarray}
Notice that for any polynomial $B$ and for $0\le t\le l-1$,
$\left(\frac{\partial}{\partial x_{1}}\right)^{t}B^{l}$ is a
multiple of $B$. Using (\ref{elong-exp}) we have
\begin{eqnarray*}
& &P Q R \prod_{i=1}^{r}(x_{2}-y_{i})^{k}
\prod_{j=1}^{s}(x_{2}-z_{j})^{k}
(\sum_{t=0}^{l-1}\frac{(-1)^{t}}{t!}x_{0}^{n+t-k}
\left(\frac{\partial}{\partial x_{1}}\right)^{t}
\left(\prod_{i=1}^{r}(x_{1}-y_{i})^{k}
\prod_{j=1}^{s}(x_{1}-z_{j})^{k}\right)^{l}\nonumber\\
&&\cdot (x_{1}-x_{2})^{k}\psi_{1}(y_{1})\cdots \psi_{r}(y_{r})
a(x_{1})b(x_{2})\phi_{1}(z_{1})\cdots
\phi_{s}(z_{s}))|_{x_{1}=x_{2}+x_{0}}\nonumber\\
&\in& (\Hom (W,W((y_{1}^{\frac{1}{N}},\dots,
y_{r}^{\frac{1}{N}},x_{2}^{^{\frac{1}{N}}},z_{1}^{\frac{1}{N}},\dots,z_{s}^{\frac{1}{N}})))
((x_{0})).
\end{eqnarray*}
Then
\begin{eqnarray}
& &P Q R\prod_{i=1}^{r}(x_{2}-y_{i})^{k(l+1)}\nonumber\\
&&\cdot \prod_{j=1}^{s}(x_{2}-z_{j})^{k(l+1)}
\psi_{1}(y_{1})\cdots\psi_{r}(y_{r})
(a(x_{2})_{n}b(x_{2}))\phi_{1}(z_{1})\cdots\phi_{s}(z_{s})\nonumber\\
&\in& \Hom (W,W((y_{1}^{\frac{1}{N}},\dots,
y_{r}^{\frac{1}{N}},x_{2}^{^{\frac{1}{N}}},z_{1}^{\frac{1}{N}},\dots,z_{s}^{\frac{1}{N}}))).
\end{eqnarray}
This proves that the sequence
$(\psi_{1}(x),\dots,\psi_{r}(x),a(x)_{n}b(x), \phi_{1}(x),
\dots,\phi_{s}(x))$ is compatible.
\end{proof}

Now we have:

\bt{tmaximal} Let $V$ be a maximal $\Z_{N}$-graded compatible
subspace of $\E(W,N)$. Then $V$ contains $1_{W}$ and is closed.
Furthermore, $(V,Y_{\cal{E}},1_{W})$ carries the structure of a
nonlocal vertex algebra with $\theta$ as an automorphism of period
$N$ and $W$ is a $\theta$-twisted $V$-module with
$Y_{W}(a(x),x_{0})=a(x_{0})$ for $a(x)\in V$. \et

\begin{proof}
Clearly, the linear span of $V$ and $1_{W}$ is still $\Z_{N}$-graded
and compatible. With $V$ being maximal we must have $1_{W}\in V$.
Let $a(x),b(x)\in V$ and $n\in \Z$. In view of Proposition
\ref{pgeneratingcomplicatedone}, any (ordered) sequence in $V\cup \{
a(x)_{n}b(x)\}$ with one appearance of $a(x)_{n}b(x)$ is compatible.
By induction on the number of appearance of $a(x)_{n}b(x)$, using
Proposition \ref{pgeneratingcomplicatedone} again, we see that any
(ordered) sequence in $V\cup \{ a(x)_{n}b(x)\}$ with any (finite)
number of appearance of $a(x)_{n}b(x)$ is compatible. Thus $V+\C
a(x)_{n}b(x)$ is compatible. Since $V$ is maximal, we must have
$V+\C a(x)_{n}b(x)= V$,  which implies that $a(x)_{n}b(x)\in V$.
This proves that $V$ is closed. By Theorem \ref{tclosed}
$(V,Y_{\cal{E}},1_{W})$ carries the structure of a nonlocal vertex
algebra with with $\theta$ as an automorphism of period $N$ and $W$
is a $\theta$-twisted $V$-module.
\end{proof}

Furthermore, we have (cf. \cite{li-g1}, \cite{li-qva1}):

\bt{tgeneratingthem} Let $U$ be any compatible subset of $\E(W,N)$,
which linearly spans a graded subspace. There exists a unique
smallest closed $\Z_{N}$-graded compatible subspace denoted by
$\<U\>$ which contains $U$ and $1_{W}$. Furthermore,
$(\<U\>,Y_{\E},1_{W})$ carries the structure of a nonlocal vertex
algebra with $\theta$ as an automorphism of period $N$ and $W$ is a
$\theta$-twisted module with $Y_{W}(a(x),x_{0})=a(x_{0})$ for
$a(x)\in\<U\>$. \et

\begin{proof} From assumption, $U$ and $1_{W}$ linearly span
a $\Z_{N}$-graded compatible subspace of $\E(W,N)$. In view of
Zorn's lemma, there exists a maximal $\Z_{N}$-graded compatible
subspace $V$ of $\E(W,N)$, containing $U$ and $1_{W}$. By Theorem
\ref{tmaximal}, $V$ is closed, $(V,Y_{\cal{E}},1_{W})$ carries the
structure of a nonlocal vertex algebra with $\theta$ as an
automorphism  of period $N$, and $W$ is a $\theta$-twisted
$V$-module. Then the nonlocal vertex subalgebra $\<U\>$ of $V$,
generated by $U$, is the unique smallest closed $\Z_{N}$-graded
compatible subspace that contains $U$ and $1_{W}$. The rest is
clear.
\end{proof}

Next, we are going to show that compatible subsets of a certain type
generate weak quantum vertex algebras.

\bd{dquasi-s-local} {\em A subset $U$ of $\E(W,N)$ is said to be
{\em $\S$-local} if for any $a(x),b(x)\in U$, there exist
$f_{i}(x)\in\C((x))$, $a_{i}(x),b_{i}(x)\in U$ for $i=1,\ldots,r$
(finite), and a nonnegative integer $k$ such that
\begin{eqnarray}\label{esr-local-simple}
& &(x_{1}-x_{2})^{k}a(x_{1})b(x_{2}) =(x_{1}-x_{2})^{k}
\sum_{i=1}^{r}f_{i}(x_{2}-x_{1}) a_{i}(x_{2})b_{i}(x_{1}).
\end{eqnarray}
}\ed

\bl{Slocal-property} Every $\S$-local subset $U$ of $\E(W,N)$ is
compatible. \el

\begin{proof} It follows from the same proof of
Lemma 3.2 in \cite{li-qva1}.
\end{proof}

The following is the main result of this section:

\bt{tSlocality-key} Let $U$ be any $\S$-local subset of $\E(W,N)$,
which linearly spans a graded subspace. Then the nonlocal vertex
algebra $\<U\>$ generated by $U$ is a weak quantum vertex algebra on
which $\theta$ acts as an automorphism of period $N$, and $W$ is a
faithful $\theta$-twisted $\<U\>$-module with
$Y_{W}(a(x),x_{0})=a(x_{0})$ for $a(x)\in \<U\>$. \et

\begin{proof} {}By Lemma \ref{Slocal-property} and  Theorem \ref{tgeneratingthem},
$\<U\>$ is a nonlocal vertex algebra with an automorphism $\theta$
of period $N$, and $W$ is a faithful $\theta$-twisted $\<U\>$-module
with $Y_{W}(a(x),x_{0})=a(x_{0})$ for $a(x)\in \<U\>$. As $U$ is
$\S$-local, it follows from Propositions \ref{pslocal-twisted} and
\ref{pwqva-generating} that $\<U\>$ is a weak quantum vertex
algebra.
\end{proof}

\section{Twisted modules for quantum vertex algebra $V_{\Q}$}
In this section we shall use the general construction we have
established in Section 3 to construct twisted modules for quantum
vertex algebras $V_{\Q}$ which were constructed in \cite{kl}.

First we recall the quantum vertex algebra $V_{\Q}$ from \cite{kl}.
Let $r$ be a positive integer and let $\Q=(q_{ij})_{i,j=1}^{r}$ be a
square matrix of complex numbers such that
$$q_{ij}q_{ji}=1 \;\;\mbox{for}\; 1\leq i,j\leq r.$$
Define $\A_{\Q}$ to be the associative algebra (over $\C$) with
identity, with generators $$X_{i,n}, Y_{i,n} \;\;\mbox{for} \;1\le
i\le r,\;n\in\Z,$$ subject to relations
\begin{eqnarray}\label{eaq-relations}
& & X_{i,m}X_{j,n}=q_{ij}X_{j,n}X_{i,m},\;\;\;
Y_{i,m}Y_{j,n}=q_{ij}Y_{j,n}Y_{i,m},\nonumber\\
&&X_{i,m}Y_{j,n}-q_{ji}Y_{j,n}X_{i,m}=\delta_{i,j}\delta_{m+n+1,0}
\end{eqnarray}
for $1\le i,j\le r,\;m,n\in\Z.$ For $1\le i\le r$, form the
generating functions
$$X_{i}(z)=\sum_{n\in \Z}X_{i,n}z^{-n-1},\ \ \
Y_{i}(z)=\sum_{n\in \Z}Y_{i,n}z^{-n-1}\in \A_{\Q}[[z,z^{-1}]].$$ The
defining relations in (\ref{eaq-relations}) now read as
\begin{eqnarray}
& & X_{i}(z_{1})X_{j}(z_{2})=q_{ij}X_{j}(z_{2})X_{i}(z_{1}),\;\;\;
Y_{i}(z_{1})Y_{j}(z_{2})=q_{ij}Y_{j}(z_{2})Y_{i}(z_{1}),\nonumber\\
&&X_{i}(z_{1})Y_{j}(z_{2})-q_{ji}Y_{j}(z_{2})X_{i}(z_{1})
=\delta_{i,j}z_{2}^{-1}\delta\left(\frac{z_{1}}{z_{2}}\right).
\end{eqnarray}

A vector $w$ in an $\A_{\Q}$-module is called a {\em vacuum vector}
if $X_{i,n}w=Y_{i,n}w=0$ for $1\le i\le r,\; n\ge 0$, and an
$\A_{\Q}$-module $W$ equipped with a vacuum vector which generates
$W$ is called a {\em vacuum $\A_{\Q}$-module}.

Denote by $J_{\bf Q}$ the left ideal of $\A_{\Q}$ generated by the
elements
$$X_{i,n},\ Y_{i,n}\ \ (\;1\le i\le r,\;n\ge 0).$$
Furthermore, set
$$V_{\Q}=\A_{\Q}/J_{\Q},$$
a left $\A_{\Q}$-module, and set
$${\bf{1}}=1+J_{\Q}\in V_{\Q}.$$ Then
$\bf{1}$ is a vacuum vector and $V_{\Q}$ with ${\bf 1}$ is a vacuum
$\A_{\Q}$-module. Moreover, the vacuum $\A_{\Q}$-module $V_{\Q}$ is
universal in the obvious sense. For $1\le i\le r$, set
\begin{eqnarray}
u^{(i)}=X_{i,-1}{\bf 1},\ \ \ \ v^{(i)}=Y_{i,-1}{\bf 1}\in V_{\Q}.
\end{eqnarray}
It was proved in \cite{kl} (cf. \cite{li-qva2}) that there exists a
unique quantum vertex algebra structure on $V_{\Q}$ with ${\bf 1}$
as the vacuum vector and with
$$Y(u^{(i)},z)=X_{i}(z),\ \ \ \ Y(v^{(i)},z)=Y_{i}(z) \ \ \ \
\mbox{ for }1\le i\le r.$$  $V_{\Q}$ as a nonlocal vertex algebra is
generated by the vectors $u^{(i)}, v^{(i)}$ for $1\le i\le r$.
Furthermore, $V_{\bf Q}$ is a conformal quantum vertex algebra of
central charge $-(q_{11}+\cdots +q_{ll})$ with conformal vector
\begin{eqnarray}
\omega=\frac{1}{2}\sum_{i=1}^{r}(v^{(i)}_{-2}u^{(i)}
-q_{ii}u^{(i)}_{-2}v^{(i)})
\end{eqnarray}
and $V_{\bf Q}$ is $\frac{1}{2}\N$-graded by $L(0)$-weights, where
$(V_{\bf Q})_{(0)}=\C{\bf 1}$ and
$$(V_{\bf Q})_{(1/2)}={\rm span}\{u^{(i)},v^{(i)}\;|\; 1\le
i\le r\}.$$

Let $N$ be a positive integer. By observing the defining relations
of $\A_{\Q}$, it is readily to see that $\A_{\Q}$ admits a unique
automorphism $\theta_{N}$ such that
\begin{eqnarray}
& & \theta_{N}(X_{i,n})=\omega_{N} X_{i,n},\;\;\;
\theta_{N}(Y_{i,n})=\omega_{N}^{-1}Y_{i,n} \ \ \ \ \mbox{ for }1\le
i\le r,\; n\in\Z.
\end{eqnarray}
As $\theta_{N} (J_{\Q})=J_{\Q}$, we see that $\theta_{N}$ gives rise
to a linear automorphism, which is also denoted by $\theta_{N}$, of
$V_{\Q}$. We have
$$ \theta_{N}({\bf 1})={\bf 1},\ \ \ \
\theta_{N} (u^{(i)}_{n}w)=\omega_{N}u^{(i)}_{n}\theta_{N}(w), \ \ \
\theta_{N} (v^{(i)}_{n}w)=\omega_{N}^{-1}v^{(i)}_{n}\theta_{N}(w)$$
for $1\le i\le r,\; n\in \Z,\; w\in V_{\Q}$, as
$u^{(i)}_{n}=X_{i,n}$ and $v^{(i)}_{n}=Y_{i,n}$.
 Since $V_{\Q}$ as a nonlocal vertex
algebra is generated by the vectors $u^{(i)}, v^{(i)}$ for $1\le
i\le r$, it follows that $\theta_{N}$ is an automorphism of the
quantum vertex algebra $V_{\Q}$, where
\begin{eqnarray}\label{automorphism}
& & \theta_{N}(u^{(i)})=\omega_{N}u^{(i)},\ \ \ \ \
\theta_{N}(v^{(i)})=\omega_{N}^{-1}v^{(i)} \ \ \ \mbox{ for }1\le
i\le r.
\end{eqnarray}
{}From now on, we fix this order-$N$ automorphism $\theta_{N}$ of
$V_{\Q}$.

Next, we study $\theta_{N}$-twisted $V_{\Q}$-modules. To describe
$\theta_{N}$-twisted $V_{\Q}$-modules we shall need a
$\theta_{N}$-twisted analogue of $\A_{\Q}$.

\bd{associ-alg} {\em Define $\A_{\Q}[\theta_{N}]$ to be the
associative algebra with identity, with generators
$$X^{\theta}_{i,m},\; Y^{\theta}_{i,n} \ \ \ \mbox{for } \;1\le i\le r,\; m\in
\frac{1}{N}+\Z,\; n\in -\frac{1}{N}+\Z,$$ subject to relations
\begin{eqnarray}\label{eaqtheta-relations}
& &X^{\theta}_{i,m}X^{\theta}_{j,n}
=q_{ij}X^{\theta}_{j,n}X^{\theta}_{i,m},\ \ \ \
Y^{\theta}_{i,m}Y^{\theta}_{j,n}=q_{ij}Y^{\theta}_{j,n}Y^{\theta}_{i,m},\nonumber\\
&&X^{\theta}_{i,m}Y^{\theta}_{j,n}-q_{ji}Y^{\theta}_{j,n}X^{\theta}_{i,m}
=\delta_{i,j}\delta_{m+n+1,0}
\end{eqnarray}
for $1\le i,j\le r,\;m\in\frac{1}{N}+\Z,\;n\in -\frac{1}{N}+\Z.$}
\ed

For $1\le i\le r$, form the generating functions
\begin{eqnarray*}
&
&X_{i}^{\theta}(x)=\sum_{m\in\frac{1}{N}+\Z}X^{\theta}_{i,m}x^{-m-1},\;\;\;
Y_{i}^{\theta}(x)=\sum_{n\in
-\frac{1}{N}+\Z}Y^{\theta}_{i,n}x^{-n-1},
\end{eqnarray*}
which are elements of $\A_{\Q}[\theta_{N}][[x^{1/N},x^{-1/N}]]$. The
relations (\ref{eaqtheta-relations}) now read as
\begin{eqnarray} & &
X_{i}^{\theta}(x_{1})X_{j}^{\theta}(x_{2})
\label{x-identity}=q_{ij}X_{j}^{\theta}(x_{2})X_{i}^{\theta}(x_{1}),\
\ \ \ \ \label{y-identity}Y_{i}^{\theta}(x_{1})Y_{j}^{\theta}(x_{2})
=q_{ji}Y_{j}^{\theta}(x_{2})Y_{i}^{\theta}(x_{1}),\nonumber\\
&&X_{i}^{\theta}(x_{1})Y_{j}^{\theta}(x_{2})-q_{ji}Y_{j}^{\theta}(x_{2})X_{i}^{\theta}(x_{1})
\label{xy-identity}=\delta_{i,j}x_{1}^{-1}\delta\left(\frac{x_{2}}{x_{1}}\right)
\left(\frac{x_{2}}{x_{1}}\right)^{\frac{1}{N}}.
\end{eqnarray}

We say that an $\A_{\Q}[\theta_{N}]$-module $W$ is {\em restricted}
if for any $w\in W$, $1\le i \le r$, $X_{i,m}w=Y_{i,n}w=0$ for $m\in
\frac{1}{N}+\Z,\; n\in -\frac{1}{N}+\Z$ sufficiently large.

\bt{tqva-construction} Let $W$ be any restricted
$\A_{\Q}[\theta_{N}]$-module. Then there exists a unique structure
of a $\theta_{N}$-twisted $V_{\Q}$-module on $W$ with
$$Y_{W}(u^{(i)},z)=X^{\theta}_{i}(z),\ \ Y_{W}(v^{(i)},z)=Y^{\theta}_{i}(z)$$
for $1\le i\le r$. On the other hand, for any $\theta_{N}$-twisted
$V_{\Q}$-module $(W,Y_{W})$, $W$ is a restricted
$\A_{\Q}[\theta_{N}]$-module with
$$X^{\theta}_{i}(z)=Y_{W}(u^{(i)},z),\ \ Y^{\theta}_{i}(z)=Y_{W}(v^{(i)},z)
\ \ \ \ \mbox{ for }1\le i\le r.$$
 \et

\begin{proof} The uniqueness is clear as $V_{\bf Q}$ is generated by
$u^{(i)},v^{(i)}\;(1\le i\le r)$. Set
$$U_{W}=\{X_{i}^{\theta}(x), Y_{i}^{\theta}(x)\;|\;1\le i\le r\}\subset \E(W,N).$$
It follows from (\ref{xy-identity}) that $U_{W}$ is an $\S$-local
subset of $\E(W,N)$, consisting of homogeneous elements. By Theorem
\ref {tSlocality-key}, $U_{W}$ generates a weak quantum vertex
algebra $\<U_{W}\>$ on which $\theta_{N}$ acts as an automorphism,
and $W$ is a faithful $\theta_{N}$-twisted $\<U_{W}\>$-module. With
the relations (\ref{xy-identity}), it follows from  Proposition
\ref{q-communication} that $\<U_{W}\>$ is an $\A_{\Q}$-module with
$X_{i,n}$ and $Y_{i,n}$ acting as $X_{i}^{\theta}(x)_{n}$ and
$Y_{i}^{\theta}(x)_{n}$ for $1\le i\le r,\; n\in \Z$, respectively.
Since
 $\<U_{W}\>$ as a nonlocal vertex algebra is generated by
$X_{i}^{\theta}(x), Y_{i}^{\theta}$ for $1\le i\le r$, it follows
from Lemma 2.3 that $\<U_{W}\>$ is a vacuum $\A_{\Q}$-module with a
vacuum vector $1_{W}$. Consequently, there exists an
$\A_{\Q}$-module homomorphism $\psi$ from $V_{\Q}$ to $\<U_{W}\>$,
sending ${\bf{1}}$ to $1_{W}$. Again, since $V_{\Q}$ is generated by
$u^{(i)}, v^{(i)}$ for $1\le i\le r$, it follows that $\psi$ is a
homomorphism of nonlocal vertex algebras. Consequently, $W$ has a
$\theta_{N}$-twisted $V_{\Q}$-module structure as desired.

Now, let $(W,Y_{W})$ be a $\theta_{N}$-twisted $V_{\Q}$-module. By
Proposition \ref{q-communication} we have
\begin{eqnarray*} & &
Y_{W}(u^{(i)},x_{1})Y_{W}(u^{(j)},x_{2})
=q_{ij}Y_{W}(u^{(j)},x_{2})Y_{W}(u^{(i)},x_{1}),\\
&&Y_{W}(v^{(i)},x_{1})Y_{W}(v^{(j)},x_{2})
=q_{ji}Y_{W}(v^{(j)},x_{2})Y_{W}(v^{(i)},x_{1}),\nonumber\\
&&Y_{W}(u^{(i)},x_{1})Y_{W}(v^{(j)},x_{2})-q_{ji}Y_{W}(v^{(j)},x_{2})Y_{W}(u^{(i)},x_{1})
=\delta_{i,j}x_{1}^{-1}\delta\left(\frac{x_{2}}{x_{1}}\right)
\left(\frac{x_{2}}{x_{1}}\right)^{\frac{1}{N}}.
\end{eqnarray*}
Thus $W$ is a restricted $\A_{\Q}[\theta_{N}]$-module with
$$X^{\theta}_{i}(z)=Y_{W}(u^{(i)},z),\ \ Y^{\theta}_{i}(z)=Y_{W}(v^{(i)},z)
\ \ \ \ \mbox{ for }1\le i\le r.$$ This completes the proof.
\end{proof}

\br{rvirasoro} {\em Let $W$ be any restricted
$\A_{\Q}[\theta_{N}]$-module. In view of Theorem
\ref{tqva-construction}, there exists a (unique)
$\theta_{N}$-twisted $V_{\bf Q}$-module structure $Y_{W}$ on $W$
such that
$$Y_{W}(u^{(i)},x)=X_{i}^{\theta}(x),\ \ \ \
Y_{W}(v^{(i)},x)=Y_{i}^{\theta}(x)$$ for $1\le i\le r$. Recall that
$V_{\bf Q}$ is a conformal quantum vertex algebra of central charge
$-(q_{11}+\cdots +q_{ll})$ with conformal vector
\begin{eqnarray*}
\omega=\frac{1}{2}\sum_{i=1}^{r}(v^{(i)}_{-2}u^{(i)}
-q_{ii}u^{(i)}_{-2}v^{(i)}).
\end{eqnarray*}
Clearly, $\theta_{N}(\omega)=\omega$. By Corollary \ref{cvirasoro},
$W$ is a module for the Virasoro algebra of central charge
$-(q_{11}+\cdots +q_{ll})$ with $L(m)$ for $m\in \Z$ given by
$$\sum_{m\in \Z}L(m)x^{-m-2}=Y_{W}(\omega,x)=\sum_{n\in
\Z}\omega_{n}x^{-n-1}.$$  By using the twisted Jacobi identity one
can express $L(n)$ for $n\in \Z$ in terms of
$X_{i,m+1/N}^{\theta},\; Y_{i,m-1/N}^{\theta}$ for $m\in \Z,\; 1\le
i\le r$.}\er

\end{document}